\documentclass[sort&compress]{elsarticle}
\usepackage{lineno,hyperref}
\modulolinenumbers[5]
\journal{Applied Mathematics and Computation}

\usepackage{epstopdf}
\usepackage{amsmath}
\usepackage{float}
\usepackage{color}
\usepackage{caption}
\usepackage{subcaption}

\begin{document}

\begin{frontmatter}
\title{Improved Probabilistic Principal Component Analysis for Application to Reduced Order Modeling }

\author{Indika Udagedara}
\address{Department of Mathematics, SC456, Clarkson University, P.O. Box 5815, Potsdam, NY, 13699-5815}
\ead{udagedig@clarkson.edu}

\author{Brian Helenbrook}
\address{Department of Mathematics and Department of Mechanical and Aeronautical Engineering, 357 CAMP, Clarkson University, P.O. Box 5725, Potsdam, NY, 13699-5725}
\ead{helenbrk@clarkson.edu}

\author{Aaron Luttman}
\address{National Security Technologies LLC, P.O. Box 98521, M/S NLV078, Las Vegas, NV, 89030}
\ead{luttmaab@nv.doe.gov}

\author{Jared Catenacci}
\address{National Security Technologies LLC, P.O. Box 98521, M/S NLV078, Las Vegas, NV, 89030}
\ead{CatenaJW@nv.doe.gov}

\begin{abstract}
In our previous work, a reduced order model (ROM) for a stochastic system was made, where noisy data was projected onto principal component analysis (PCA)-derived basis vectors to obtain an accurate reconstruction of the noise-free data.  That work used techniques designed for deterministic data; PCA was used for the basis function generation and $L_2$ projection was used to create the reconstructions.   In this work, probabilistic approaches are used.  The probabilistic PCA (PPCA) is used to generate the basis, which then allows  the noise in the training data to be estimated.  PPCA has also been improved so that the derived basis vectors are orthonormal and the variance of the basis expansion coefficients over the training data set can be estimated.  The standard approach assumes a unit variance for these coefficients.  Based on the results of the PPCA, model selection criteria are applied to automatically choose the dimension of the ROM.  In our previous work, a heuristic approach was used to pick the dimension.   Lastly, a new statistical approach is used for the projection step where the variance information obtained from the improved PPCA is used as a prior to improve the projection.  This gives improved accuracy over $L_2$ projection when the projected data is noisy.  In addition, the noise statistics for the projected data are not assumed to be the same as that of the training data, but are estimated in the projection process.  The entire approach gives a fully stochastic method for computing a ROM from noisy training data, determining ideal model selection, and projecting noisy test data, thus enabling accurate predictions of noise-free data  from data that is dominated by noise.
\end{abstract}

\end{frontmatter}

\section{Introduction}
In our previous work \cite{udagedara2015reduced}, we developed a reduced order modeling (ROM) technique that could be applied to problems with a stochastic component.  In that work, the application was to radiation transport problems.  The approach to generate the model was to first obtain a set of low-noise training data corresponding to different physical scenarios.  In the radiation transport problem, this was done by performing high fidelity and computationally intensive Monte-Carlo simulations to generate radiation energy spectra at different locations relative to the radiation source for different radiation source mixtures.  Principal component analysis (PCA) \cite{peason1901lines,wold1987principal,abdi2010principal,richardson2009principal, bro2014principal, jolliffe2002principal, moore1981principal} (or equivalently the proper orthogonal decomposition (POD \cite{liang2002proper,chatterjee2000introduction,volkwein2013proper,kerschen2002physical, willcox2002balanced, berkooz1993proper} ) was then applied to generate a set of basis functions that could compactly represent all of the training data.   This basis  was then used to make accurate estimates of the radiation energy spectrum given noisy data from some new scenario.  Depending on the application, the noisy data, which we call the ``trial'' data, could come from experimental measurement or a low fidelity stochastic simulation.  The ROM estimates were obtained using an $L_2$ projection of this trial data onto the space spanned by the basis functions.  In the radiation transport problem, we found that the energy spectra associated with a new mixture of materials could be calculated with four orders of magnitude reduction in the computational cost required  relative to that of a training data simulation.

Although \cite{udagedara2015reduced} did define the basic procedure for generating a ROM of a stochastic problem, there were several issues that were not addressed.  The formulation in \cite{udagedara2015reduced} did not properly account for noise in the generation of the basis functions, had no mechanism for selecting the number of basis functions in the ROM, and used $L_2$ projection of the trial data without mathematical justification.  In this work, a probabilistic approach is used to remedy these deficiencies.  To generate the basis,  the stochastic formulation of PCA, known as probabilistic principal component analysis (PPCA)  \cite{lawrence2005probabilistic, tipping1999mixtures, tipping1999probabilistic, beckmann2004probabilistic, wu2008gene, ilin2010practical}, is used.  PPCA  identifies the noise in the training data, which cannot be estimated with conventional PCA.   We have also improved PPCA to relax the standard assumption that the latent variables have unit variance.  The new formulation provides more physical insight into the eigenvalues of the PPCA and their relation to the variance of the latent variables.  This information was necessary for the new projection procedure that was developed.

The probabilistic formulation also provides a method for selecting the number of basis functions to include in the ROM. This can be done by combining PPCA with Bayesian model selection \cite{chen2008extended, toni2009approximate, wasserman2000bayesian, burnham2003model, zucchini2000introduction, beck2004model, raftery1993bayesian, raftery1995bayesian, posada2004model, johnson2004model} criterion.  Here the Bayesian information criteria (BIC) \cite{chen2008extended,burnham2004multimodel, konishi2004bayesian, claeskens2008model, wasserman2000bayesian, johnson2004model}  is used to identify the optimal number of basis functions to include in the ROM, and it is demonstrated that this approach reliably chooses the number of basis functions that can be identified given that the training data itself includes noise.

Lastly, a new approach for projecting the trial data is derived.  This approach uses prior information obtained from the PPCA of the training data to improve the projection of the trial data.  In our previous work, $L_2$ projection was used, which basically corresponds to a projection with no prior knowledge of the projection coefficients.  In the following, it is demonstrated that using the training data to define a prior for the projection coefficients leads to significantly improved results when the trial data has more noise than the training data, which is typically the case.

The paper is organized as follows.  First the framework for the reduced order model generation and the modified PPCA approach is derived.  This is then followed by a discussion of the model selection approach (BIC) and the derivation of the method for projecting the trial data.   To demonstrate the benefits of this new probabilistic formulation, the ROM procedure is applied to a simple stochastic model problem, and predictions from the ROM are compared to our previous approach and to noise-free data, which was known for the model problem.

\section{Reduced Order Modeling Formulation}\label{ROMform}
The ROM is formulated assuming the data is created by a process of the following form
\begin{equation} \label{eq1}
\vec{y} =  \sum_{j=1}^{m}w_j\vec{\phi_j} + \vec{\mu} + \vec{\epsilon}  =  \Phi \vec{w} + \vec{\mu} + \vec{\epsilon},
\end{equation}
where $\vec{y}$ is a single data realization for a particular physical scenario which could be obtained either from a numerical simulation or an experiment.  The dimension of $\vec{y}$ is $d$, which depends on the physical problem being studied.  For example, in the radiation transport problem, $d$ was the number of energy bins used to describe the radiation energy spectrum.    The functions $\vec{\phi}_{j},\ j \in [1,m]$ are basis functions that are scenario independent and the $w_{j}$'s are latent variables that vary depending on the physical scenario. The $w_{j}$'s are assumed to be random variables.  This implies that there is a probability associated with the occurrence of each physical scenario.  The mean of the latent variables is assumed to be zero such that $\vec{\mu}$ is the mean of the data over all scenarios.  $\vec{\mu}$ is thus a scenario independent constant vector.  $\vec{\epsilon}$ is a random variable that represents the noise in the process.  This noise could represent noise in the experimental measurements or in the stochastic numerical simulation approach used to generate the data.  The noise is assumed to be generated by a zero-mean Gaussian process with covariance $\sigma _{\epsilon}^{2} I $ where $I$ is the identity matrix of dimension $d$.  

Although the data is assumed to be generated by a process of the form given by \eqref{eq1}, none of the parameters of the model ($\Phi, \vec{\mu}, m$) are known.  The first step of the reduced order modeling process is to generate a set of ``training data'' that can be used to estimate these parameters.  The training data is a set $Y = \{\vec{y}_{k}\},$ for $k = 1,2, ..., n$ of realizations of the process.  The scenarios associated with these realizations are chosen randomly according to the probability density function predicting the occurrence of any given scenario.  The generation of the data can be through either numerical simulation or experiment and both are assumed to also include random noise.  

\subsection{PPCA}  \label{MPPCA}
To estimate $\Phi$ and $\vec{\mu}$ given $Y$, PPCA is used.  Our formulation of PPCA is similar to that of \cite{tipping1999mixtures,tipping1999probabilistic}.  The main difference is that in \cite{tipping1999mixtures,tipping1999probabilistic} it is assumed that the latent variables, $\vec{w}$, are uncorrelated and follow a Gaussian distribution with unit covariance.  In the following, it is also assumed that the latent variables are uncorrelated and follow a Gaussian distribution, but the covariance is not a-priori assumed to be 1.  Instead the variances, $\{\sigma^2_{w_i}\}_{i=1}^{m} $, are estimated as part of the derivation.  The derivation is similar to the original derivation provided in \cite{tipping1999mixtures,tipping1999probabilistic} so in the following, a condensed derivation is given that highlights the main differences.

Bayes' formula to estimate the unknowns, $\Phi$, $\sigma_{\epsilon}^{2}$, $\sigma_{w_i}^2$, and $\vec{\mu}$ in the model is
\begin{equation}
p\left(\Phi , \sigma_{\epsilon}^{2}, \sigma_{w_i}^2, \vec{\mu} |  Y \right)  \propto p\left (Y | \Phi, \sigma_{\epsilon}^{2}, \sigma_{w_i}^2, \vec{\mu} \right) p\left(\Phi ,\sigma_{\epsilon}^{2}, \sigma_{w_i}^2, \vec{\mu}  \right),
\end{equation}
where $p \left(\Phi, \sigma_{\epsilon}^{2}, \sigma_{w_i}^2, \vec{\mu} |  Y  \right) $ is the posterior distribution, $p\left(\Phi ,\sigma_{\epsilon}^{2}, \sigma_{w_i}^2 , \vec{\mu}  \right)$  is the prior distribution and 
$p \left(Y |  \Phi , \sigma_{\epsilon}^{2}, \sigma_{w_i}^2 , \vec{\mu} \right)$ is the likelihood distribution.
The PPCA uses a maximum likelihood estimator (MLE) \cite{fisher1922mathematical, scholz1985maximum, white1982maximum, efron1978assessing} to find the unknown parameters assuming no prior knowledge about their values i.e. the prior is assumed to be uniform.  The MLE is thus obtained by maximizing the log-likelihood function, $ \log p \left( Y | \Phi , \sigma_{\epsilon}^2, \sigma_{w_i}^2, \vec{\mu}  \right) $.

In order to obtain the likelihood function, the probability distribution of an individual realization, $\vec{y}$, conditioned on $\Phi , \sigma_{\epsilon}^2, \sigma_{w_i}^2, \vec{\mu}$, first needs to be identified. This distribution, $p(\vec{y} | \Phi , \sigma_{\epsilon}^2, \sigma_{w_i}^2, \vec{\mu})$, is called the predictive distribution and can be obtained using the following relations
\begin{equation}\label{predictive}
p(\vec{y} | \Phi , \sigma_{\epsilon}^2, \sigma_{w_i}^2, \vec{\mu}) = \int_{-\infty}^\infty p(\vec{y}, \vec{w}  | \Phi , \sigma_{\epsilon}^2, \sigma_{w_i}^2, \vec{\mu}) d\vec{w} =  \int_{-\infty}^\infty p(\vec{y} | \vec{w}, \Phi , \sigma_{\epsilon}^2, \vec{\mu}) p(\vec{w} | \sigma_{w_i}^2) d\vec{w},
\end{equation}
where the integral is a multidimensional integral over all components of the vector $\vec{w}$.  Using the assumption that the noise in \eqref{eq1} is Gaussian with zero mean, the probability distribution of $\vec{y}$ conditioned on the latent variable, $\vec{w}$, and the parameters $\sigma_{\epsilon}^2$, $\vec{\mu}$, and $\Phi$ is given by
\begin{equation}
p \left( \vec{y} |  \vec{w}, \Phi, \sigma_{\epsilon}^{2}, \vec{\mu} \right) =  \mathcal{N}\left(\Phi\vec{w} + \vec{\mu}, \sigma_{\epsilon}^2I\right)
\end{equation}
where the notation $\mathcal{N}\left(\Phi\vec{w} + \vec{\mu}, \sigma_{\epsilon}^2I\right)$ indicates a Gaussian distribution with mean $\Phi\vec{w} + \vec{\mu}$ and co-variance matrix $ \sigma_{\epsilon}^2I$.  

As mentioned above, in this analysis the only difference is that the latent variables are assumed to be uncorrelated and follow a zero-mean Gaussian distribution with a covariance of $\sigma_{w_i}^2$ instead of 1, thus $p(\vec{w} | \sigma_{w_i}^2) = {\mathcal N}(\vec{0},\Sigma^2)$ where $\Sigma^2$ is an $m \times m$ diagonal matrix with the values $\sigma^2_{w_i}$ on the diagonal.   Based on this assumption,
one can show that the predictive distribution is Gaussian of the form $\mathcal{N}\left(\vec{\mu},\Phi \Sigma^2 \Phi^{T} + \sigma_{\epsilon}^2 I\right)$, which is similar to \cite{tipping1999mixtures,tipping1999probabilistic} except for the introduction of the diagonal matrix $\Sigma^2$.

The likelihood distribution, $p \left( Y | \Phi , \sigma_{\epsilon}^2, \Sigma^2, \vec{\mu}  \right)$ for the set of training data, $Y = \{\vec{y}_{k}\}$, for $k = 1,2, ..., n$ is the product of the individual predictive distributions.  The log likelihood can be shown to be
\begin{align}\label{logLH}
L &= \log{p \left( Y  |   \Phi, \sigma_{\epsilon}^{2}, \Sigma^2, \vec{\mu} \right)} & \nonumber \\
& = - \frac{dn}{2} \log{2 \pi} - \frac{n}{2} \log |\Phi \Sigma^2 \Phi^{T} + \sigma_{\epsilon}^2 I | \nonumber \\
& - \frac{1}{2} \sum_{k=1}^{n} \left(\vec{y}_{k}-\vec{\mu}\right)^T \left(\Phi \Sigma^2 \Phi^{T} + \sigma_{\epsilon}^2 I\right)^{-1} \left( \vec{y}_{k}-\vec{\mu}\right),
\end{align}
where $|\cdot|$ denotes the determinant of a matrix.  

As the prior is uniform, the most probable values of the posterior distribution can be determined by maximizing $L$ with respect to the unknown parameters, $\Phi$, $\vec{\mu}$, $\sigma_{\epsilon}^{2}$. In all the following, the subscript $MP$ indicates these most probable values.  

The maximization process gives the following estimate for the mean $\vec{\mu}$,
\begin{equation}\label{smean}
 \vec{\mu}_{MP} = \frac{1}{n} \sum_{k=1}^{n} \vec{y}_{k}.
\end{equation}
The estimate of $\sigma_{\epsilon}^2$ is
\begin{equation}\label{sig}
\sigma_{{\epsilon}\left(MP\right)}^2 = \frac{1}{d-m} \sum_{i=m+1}^{d} \lambda_{i},
\end{equation}
where the $\lambda_{i}$'s are the eigenvalues of the data covariance matrix, \\
$S = \frac{1}{n} \sum_{k=1}^{n}\left(\vec{y}_{k}-\vec{\mu}_{MP}\right) \left( \vec{y}_{k}-\vec{\mu}_{MP}\right)^T$.  The maximum likelihood estimate for the $\sigma_{{\epsilon}}^2$ can be interpreted as the average magnitude of the eigenvalues of dimension greater than $m$. These eigenvalues can only be caused by noise, as there are only $m$ latent variables. 

Minimizing \eqref{logLH} with respect to $\Phi$, the following equation can be obtained. 
\begin{equation}\label{phi}
\Phi_{MP}{\Sigma}_{MP} = U \left(\Lambda - \sigma_{\epsilon \left(MP\right)}^2 I\right)^{\frac{1}{2}} R,
\end{equation}
where $U$ is a  $d \times m$ matrix  whose columns are given by a complete subset of (orthonormal) eigenvectors of the data covariance matrix $S$, 
$\Lambda$ is the $m \times m$ diagonal matrix consisting of the first $m$ largest eigenvalues of $S$, and $R$ is an arbitrary $m \times m$ orthonormal matrix.   This equation is almost the same as in \cite{tipping1999mixtures,tipping1999probabilistic}, except for the $\Sigma$ term.  In \cite{tipping1999mixtures,tipping1999probabilistic}, $R$ was chosen to be the identity matrix, which then determined $\Phi$.  

The disadvantage of this choice is that the column vectors of $\Phi$ then each have a magnitude determined by the diagonal matrix $\left(\Lambda - \sigma_{\epsilon \left(MP\right)}^2 I\right)^{\frac{1}{2}}$.  This scaling of the basis function is necessary to ensure that the latent variables all have a variance of unity.  In the new formulation, we can satisfy \eqref{phi} by choosing
\begin{align}\label{basis_func}
\Phi_{MP} = U
\end{align}
and the estimate for $\Sigma$ as
\begin{align} \label{var_latent}
\Sigma_{MP} = \left(\Lambda - \sigma_{\epsilon \left(MP\right)}^2 I\right)^{\frac{1}{2}}.
\end{align}
This allows us to have unit basis functions, and also correctly identifies the variance of the latent variables (as we confirm in the example problem below).  Manipulating \eqref{var_latent}, the relation can put in the following form
\begin{equation}\label{Lambda}
\Lambda = \Sigma_{MP}^2 + \sigma_{\epsilon \left(MP\right)}^2 I.
\end{equation}
That shows that the eigenvalues of the covariance matrix $S$ are the variance of latent variables summed with the variance of the measurement error.  The first $m$ eigenvalues consists of both variances, however the eigenvalues greater than $m$ are strictly due to random measurement error. 

\subsection{Model Selection} \label{dimension}

In the previous section, the model parameters were estimated assuming that the dimension of the ROM, $m$, was a known, fixed number. In our previous work, no systematic method for choosing $m$ was provided.   PPCA together with Bayesian model selection criteria can be used to predict the number of basis functions required for the ROM. There are many different Bayesian model selection criteria \cite{johnson2004model, burnham2003model, wasserman2000bayesian, stephan2009bayesian, acquah2010comparison, cavanaugh1997unifying, bozdogan2000akaike, posada2004model}.  Here the Bayesian information criterion (BIC) \cite{chen2008extended,burnham2004multimodel, konishi2004bayesian, claeskens2008model, wasserman2000bayesian, johnson2004model} is used.  BIC chooses the value of $m$ that minimizes the following function
\begin{equation} \label{BIC}
f_{BIC}(m) = -2 L_{MP}+ \left(m\left(d-1-\frac{m-1}{2}\right) +d +1 \right)  \log n,
\end{equation}
where $L_{MP}$ is the maximum value of the likelihood distribution and the term in the outer parenthesis in \eqref{BIC} is the number of estimatable parameters in the model, both of which depend on $m$.  The number of estimatable parameters arise from, $\Phi_{MP}$, $\vec{\mu}_{MP}$ and $\sigma_{\epsilon_{MP}}^2$. There are $\left(d-1\right) +\left(d-2\right) + \left(d-3\right) + ...+ \left(d-1-\left(m-1\right)\right)$ parameters in $\Phi_{MP}$. Here the $d-1$ comes from the fact that the first basis vector is required to be normalized to have magnitude 1 so when $m=1$ one can only choose $d-1$ independent variables.  Because of the requirement of orthogonality the number of free parameters in choosing a basis vector decreases by 1 for each additional basis vector.  This results in the number of free parameters in $\Phi$ being  $m(d-1-(m-1)/2)$.   The number of parameters in $\vec{\mu}_{MP}$ and $\sigma_{\epsilon_{MP}}^2$ are $d$ and 1, respectively giving the total shown in parentheses above.

$L_{MP}$ is obtained by inserting the maximum likelihood estimates of the parameters \eqref{smean}, \eqref{sig}, \eqref{basis_func}, and \eqref{var_latent}) in \eqref{logLH}. Following the simplification techniques in \cite{tipping1999mixtures} but with our maximum likelihood results, this becomes

\begin{align}\label{dim_sel}
L_{MP} &=  - \frac{dn}{2} \log \left(2\pi \right) - \frac{n}{2} \left( \sum_{j=1}^{m} \log \left( \lambda_{j} \right) + \left( d-m \right) \log \left(\frac{1}{d-m} \sum_{j=m+1}^{d} \lambda_{j} \right)  + d \right).
\end{align}
To find the most probable value for $m$, $f_{BIC}(m)$ is calculated as a function of $m$,  $m \in [1,d]$,  and the value of $m$ that minimizes the function is chosen.

\subsection{Bayesian Projection with Gaussian Prior}\label{projection}
The above sections determined the model parameters, $\Phi$, $\vec{\mu}$ and $m$.  In this section, a latent variable vector $\vec{w}$ is estimated given a ``trial'' data vector $\vec{y}$ that is obtained from a new scenario drawn from the distribution of scenario probabilities.  $\vec{y}$ also includes noise, $\epsilon_T$, which is drawn from a zero-mean Gaussian distribution $\mathcal{N}(0,\sigma^2_{\epsilon_T})$ where the magnitude of this noise, $\sigma^2_{\epsilon_T}$, is assumed to be different (typically larger) than that of the training data.   In our previous work $L_2$ projection of the trial data was used to estimate the latent variables and no estimate was given for $\sigma^2_{\epsilon_T}$.   Here the latent variables are estimated using Bayesian parameter estimation with a Gaussian prior.  The estimate of $\Sigma_{MP}^2$ obtained from the training data is used as the covariance of the prior on the latent variables.

Assuming the model given by \eqref{eq1} holds for the trial data as well, the probability distribution of  $\vec{y}$ is, 
\begin{align} \label{py}
p \left( \vec{y} |  \vec{w}, \Phi, \vec{\mu}, \sigma_{\epsilon_{T}}^{2} \right) 
 &= |\sigma_{\epsilon_{T}}^{2}I |^{-1/2}  \exp\left({-\frac{\left(\vec{y} - \Phi \vec{w} -\vec{\mu} \right) ^T \left( \vec{y}- \Phi \vec{w} -\vec{\mu} \right)}{2\left(\sigma_{\epsilon_{T}}^{2}  \right)}}\right) 
\end{align}
As assumed before, the probability of $\vec{w}$ for a given scenario is Gaussian with mean zero and covariance $\Sigma^2$
\begin{align} \label{pw}
p \left(\vec{w} | \Sigma^2  \right) 
 \propto | \Sigma^2 |^{-1/2}  \exp\left({-\frac{1}{2} \vec{w}  ^T  \Sigma ^{-2}  \vec{w} }\right)
\end{align}
To estimate the probability distribution of $\vec{w}$ and $\sigma_{\epsilon_{T}}^{2}$ conditioned on the observed data, $\vec{y}$, and parameters $\vec{\mu}$, $\Sigma^2$, and $\Phi$, Bayes' theorem is used.
Applying Bayes' theorem assuming that $\Sigma^2$ and $\sigma_{\epsilon_{T}}^{2}$ are independent we have,
\begin{align}
p \left(\vec{w}, {\sigma} _{\epsilon_{T}}^{2}   | \vec{y}, \Phi, \Sigma^2,\vec{\mu}\right)  & \propto p \left(\vec{y}  |  \vec{w}, {\sigma} _{\epsilon_{T}}^{2}, \Phi, \Sigma^2, \vec{\mu} \right) {p\left(\vec{w} , {\sigma} _{\epsilon_{T}}^{2}  | \Phi, \Sigma^2,\vec{\mu} \right)}  \nonumber \\
&  \propto p \left(\vec{y}  |   \vec{w}, {\sigma} _{\epsilon_{T}}^{2}, \Phi, \vec{\mu} \right) p\left(\vec{w} | \Sigma^2  \right) p\left({\sigma} _{\epsilon_{T}}^{2}\right).
\end{align}
Assuming  a uniform prior distribution for $\sigma _{\epsilon_{T}}^{2}$, the log posterior can then be obtained from \eqref{py} and \eqref{pw} as
\begin{align}\label{logpost}
  \log p \left(\vec{w}, {\sigma} _{\epsilon_{T}}^{2}   | \vec{y}, \Phi, \Sigma^2,\vec{\mu}\right) & \nonumber \propto 
     \frac{1}{2} \log |\sigma_{\epsilon_{T}}^{2} I | + \frac{\left(\vec{y} - \Phi \vec{w} - \vec{\mu} \right)^{T} \left(\vec{y}- \Phi \vec{w} - \vec{\mu}\right)}{2\sigma_{\epsilon_{T}}^{2}}  \\
  & \qquad{}+  \frac{1}{2} \log |\Sigma^2 | + \frac{1}{2} \vec{w} ^T \Sigma^{-2}  \vec{w}
\end{align}
Setting the derivatives of $\eqref{logpost}$ with respect to $\vec{w}$ and $\sigma_{\epsilon_{T}}^{2}$ to zero to find the maximum gives
\begin{align}\label{MAP3}
\vec{w}_{{MP}}& = \left( I +  \sigma_{\epsilon_{T_{MP}}}^{2} \Sigma^{-2} \right) ^{-1} \Phi^{T}   \left( \vec{y}- \vec{\mu}  \right) 
\end{align}
and 
\begin{align}\label{eq_sigma}
\sigma_{\epsilon_{{T_{MP}}}}^{2} = \frac{1}{d} \left(\vec{y} - \Phi \vec{w}_{MP} - \vec{\mu} \right)^{T} \left(\vec{y}- \Phi \vec{w}_{MP} - \vec{\mu} \right).
\end{align}
Equation $\eqref{MAP3}$ and $\eqref{eq_sigma}$ are a system of non-linear equations with unknowns $\vec{w}_{{MP}}$ and  $\sigma_{\epsilon_{T_{MP}}}^{2}$.   The pair of equations can be solved using a fixed point iteration where it is first assumed that $\sigma^2_{\epsilon_{T_{MP}}}$ is zero, then \eqref{MAP3} is used to calculate $\vec{w}_{MP}$. Equation \eqref{eq_sigma}, which can be interpreted as a calculation of the noise in the data assuming the true data is given by $\Phi \vec{w}_{MP} + \vec{\mu}$, can then be used to calculate $\sigma^2_{\epsilon_{T_{MP}}}$. 
This new value of $\sigma^2_{\epsilon_{T_{MP}}}$ is then used in $\eqref{MAP3}$ and the process is  repeated until $\eqref{MAP3}$ and $\eqref{eq_sigma}$ are satisfied to a specified tolerance.

If $\sigma_{\epsilon_{T_{MP}}}^{2}$ is small relative to the values of $\sigma_{w_i}^2$, then the matrix in parentheses in \eqref{MAP3} is essentially the identity matrix and the $L_2$ projection result, $\vec{w} =  \Phi^{T}   \left( \vec{y}- \vec{\mu}  \right)$, is recovered.  This is the approach that was used in our previous work and is also the result that would be obtained assuming a uniform prior on $\vec{w}$ instead of a Gaussian prior.   However when the noise in the data is large relative to the variance of a latent variable, i.e. $\sigma_{\epsilon_{T_{MP}}}^{2}/\sigma_{w_i}^2$ is large, the Gaussian-prior projection reduces the magnitude of the $L_2$ projection value of $w_i$ to account for the fact that the noise in the data  is causing an estimation for $w_i$ that is larger than the expected variation of that latent variable.  In the model problem, we show that this significantly improves the accuracy of the projection for these conditions.

\section{Demonstration}\label{Model prob}
This section illustrates the above discussed ROM techniques for a model problem where the data is generated using a model of the form given by \eqref{eq1} i.e. the data is generated as a linear combination of a finite number of basis functions and latent variables with a mean vector and added random noise.  As all of the parameters of the data generation are known, the ROM process can be validated by comparing the estimated parameters to those used to generate the data.   

The $m$ basis functions used to generate the data are discrete sine waves given by
\[
\vec{\phi_j} =  \frac{\sin\left(j \pi \vec{x} \right)}{|| \sin\left( j \pi \vec{x} \right) ||} \quad \text{for} \quad j \in [1,m],
\]
where $\vec{x} \in {\mathcal R}^d$ is a vector of $d=100$ uniformly spaced points from the domain $[0,1]$ including endpoints.  In the above, the norm $|| . ||$ is the Euclidean vector norm such that $\vec{\phi_j}\cdot \vec{\phi_j} = 1$.   With this normalization, the peak value of the basis functions is $1/\sqrt(d/2) \approx 0.14$.

The latent variables, $w_j$, were sampled from Gaussian distributions with variances of $\sigma^2_{w_j}$.  The values of $\sigma_{w_j}^2$ were given by
\[
\sigma^2_{w_j} = \frac{1}{2^{j-1}} \quad \text{for} \quad j \in [1,m].
\]
Unless stated otherwise, $m = 10$ basis functions and latent variables were used to create the data.   The mean, $\vec{\mu}$, was a vector of ones.  

\begin{figure}[H]\vspace{1cm}
\centering
(a)\includegraphics[width=5.5cm]{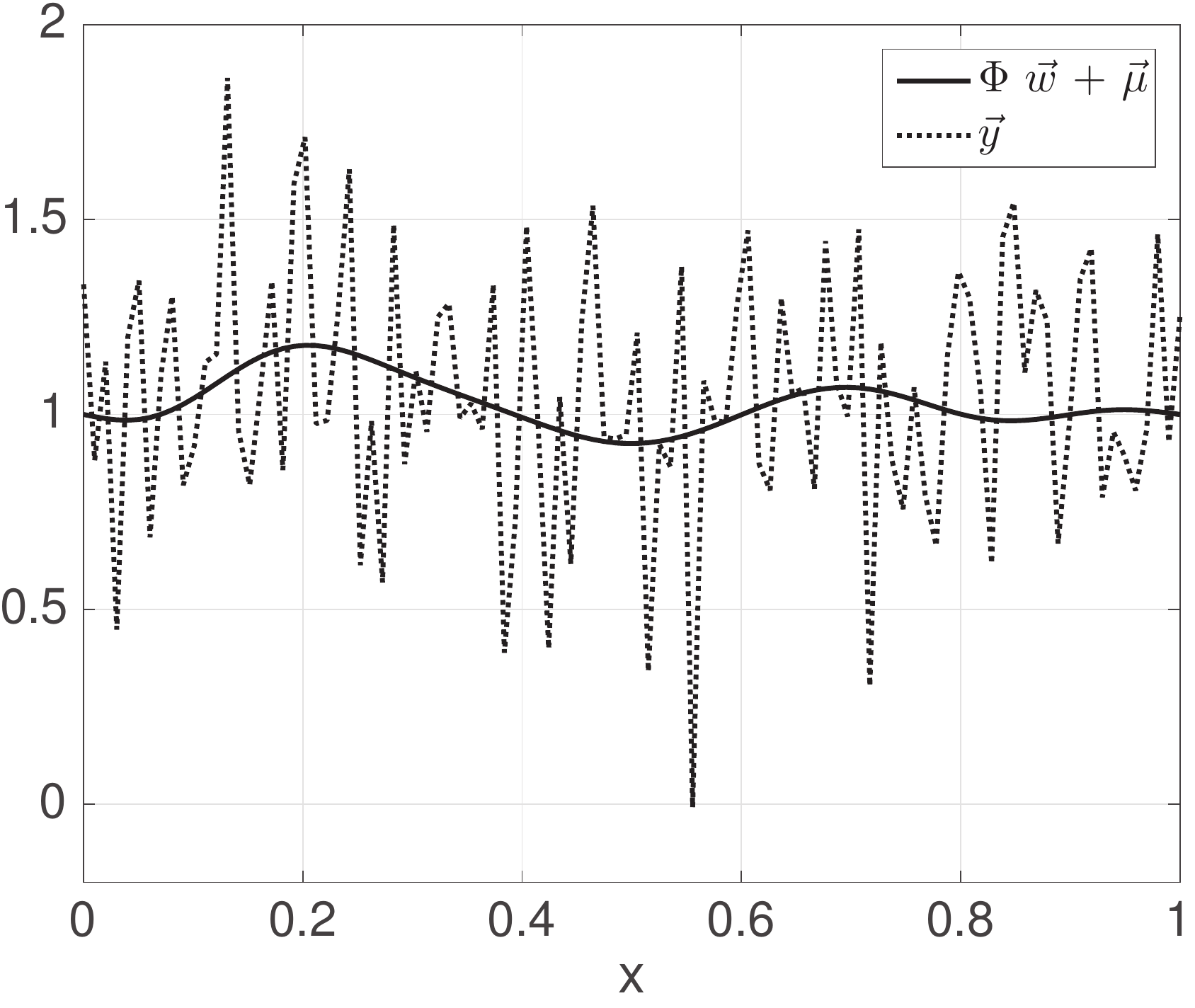}
(b)\includegraphics[width=5.5cm]{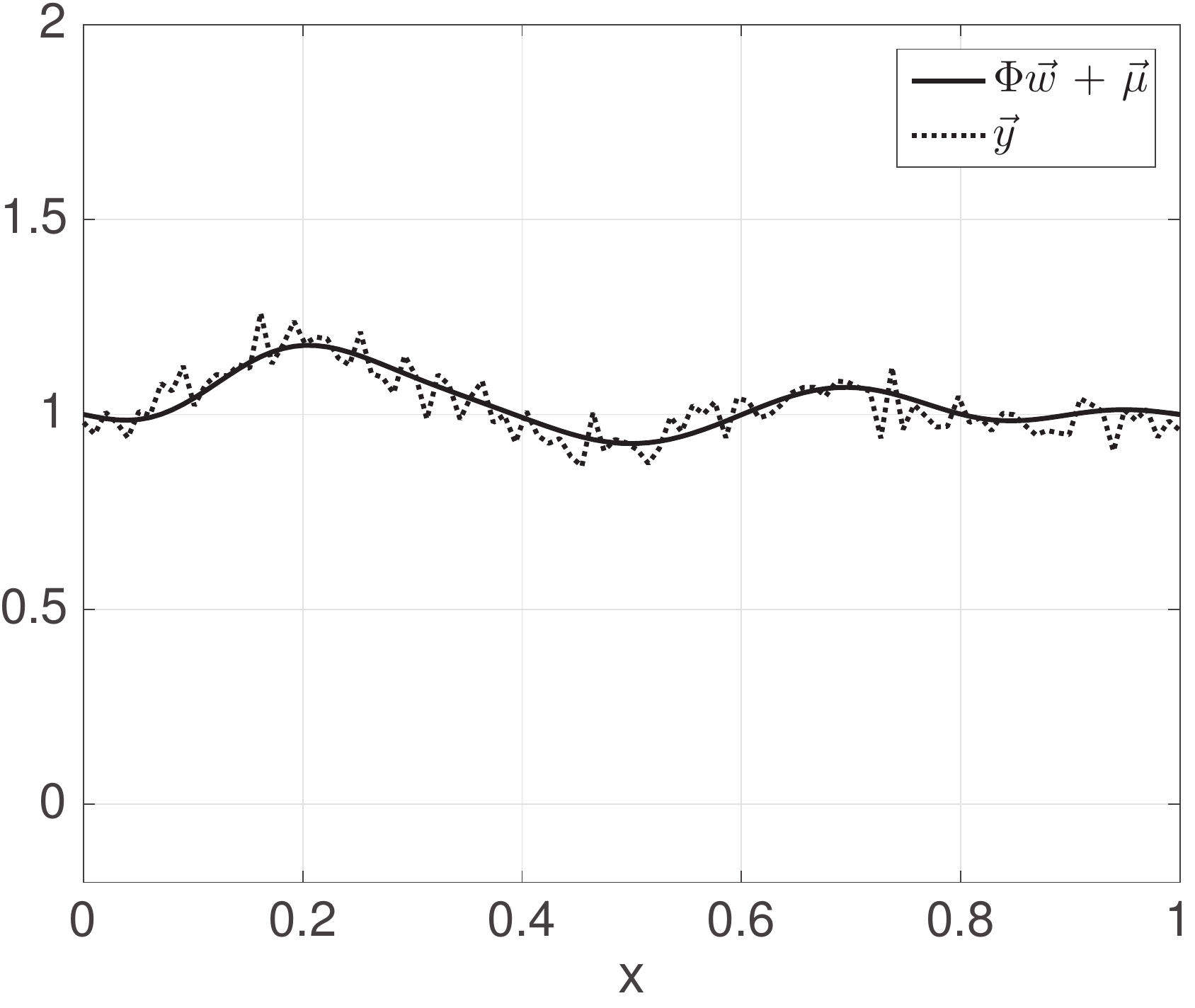}
\caption{A single realization from the training data (a) $\sigma_{\epsilon}^{2}$ = 1/10  (b) $\sigma_{\epsilon}^{2}$ = 1/400.}
\label{data_0_1}
\end{figure}

The noise vector, $\vec{\epsilon}$ was also sampled from a Gaussian distribution, either $\mathcal{N} \left( \vec{0}, \sigma_{\epsilon}^{2} I \right)$ for the training data or $\mathcal{N} \left( \vec{0}, \sigma_{\epsilon_T}^{2} I \right)$ for the trial data.   The data generation process was repeated $n = 10000$ times to generate the training data set $Y$.   To investigate the effect of noise in the training data, two training data sets were studied, one with $\sigma_{\epsilon}^{2} = 1/10$ and the other with $\sigma_{\epsilon}^{2} = 1/400$.   Fig.~\ref{data_0_1} shows a typical realization of a data vector from the two training data sets. From fig.~\ref{data_0_1}(a) it can be seen that the case of $\sigma_{\epsilon}^{2} = 1/10$  the noise in the data is significant with a magnitude on the order $\pm \sqrt{1/10} \approx 0.32$. Similarly in fig.~\ref{data_0_1}(b), the noise magnitude is $\pm \sqrt{1/400} = 0.05$.

\vspace{1cm}

\begin{figure}[H]\vspace{1cm}
 \centering
\includegraphics[width=8cm]{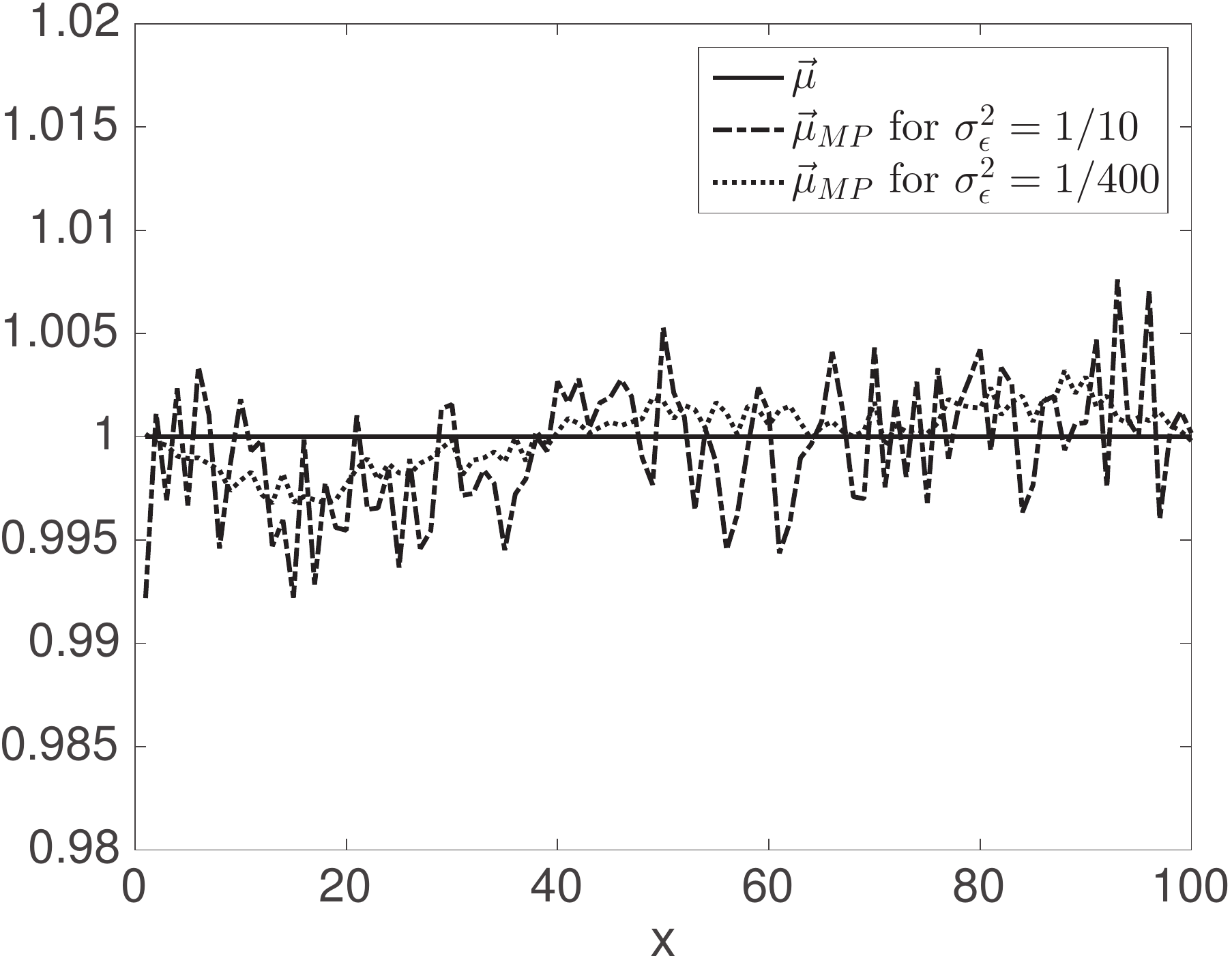}
\caption{Mean distribution for the $\sigma_{\epsilon}^2 = 1/10$ and $\sigma_{\epsilon}^2 = 1/400$ training data.}
\label{mean}
\end{figure}

Averaging over all of the data vectors of the data set $Y$ determines the most probable mean vector $\vec{\mu}_{MP}$.  This is shown in Fig.~\ref{mean} for both data sets.  The fluctuations in the mean should scale as $\left( \sum_i \sigma_{w_i} +\sigma_\epsilon\right) / \sqrt{n}$ which is equal to 0.0088 and 0.0061 for $\sigma_\epsilon^2 = 1/10$ and $1/400$ respectively.  This is in good agreement with what is observed in Fig.~\ref{mean}.   The dominant source of error in determining the mean is not the noise in the data, but rather determining the average outcome of the scenarios.  For this reason, the fluctuations in both estimated means are of similar magnitude.
 \vspace{1cm}

\subsection{PPCA}
\begin{figure}[H]\vspace{1cm}
 \centering
\includegraphics[width=8cm]{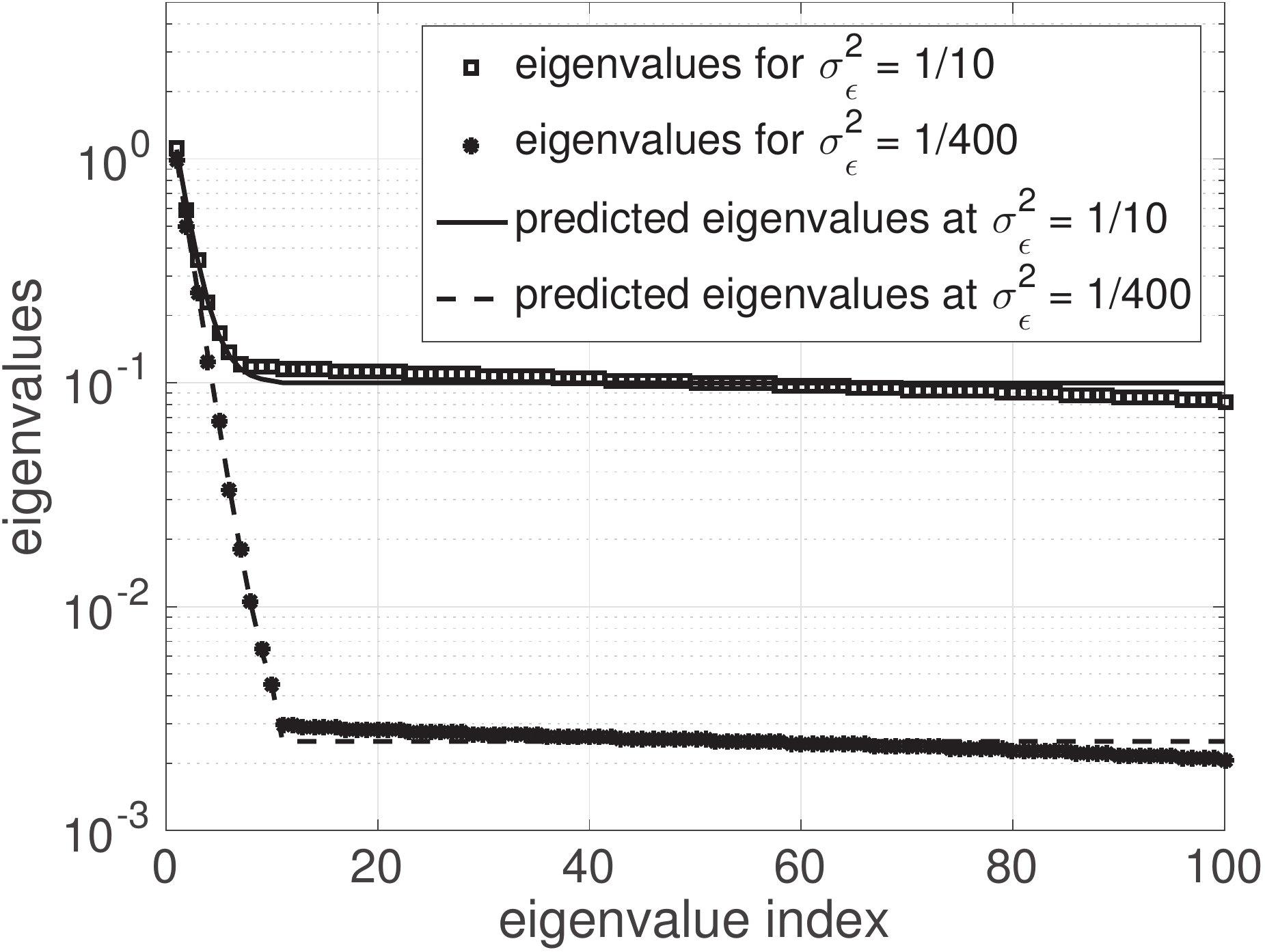}
\caption{Eigenvalue spectra of the data covariance matrices created using the $\sigma_{\epsilon}^2 = 1/10$ and $\sigma_{\epsilon}^2 = 1/400$ training data. Predicted spectra based on \eqref{Lambda} are also shown.}
\label{eig} 
\end{figure}

According to \eqref{basis_func} in Section \ref{MPPCA}, the estimate for the most probable basis functions are the dominant subset of $m$ eigenvectors of the data covariance matrix, $S = \frac{1}{n} \sum_{k=1}^{n}\left(\vec{y}_{k} -\vec{\mu}_{MP}\right) \left( \vec{y}_{k}-\vec{\mu}_{MP}\right)^T$. Fig.~\ref{eig} shows the eigenvalue spectrum for the two sets of training data with $\sigma_\epsilon^2 = 1/10$ and $1/400$.  Based on \eqref{Lambda}, the eigenvalues are expected to decay like $1/2^{i-1}$ because of the $\sigma^2_{w_i}$ term and then plateau at a value of $\sigma_\epsilon^2$.  The predicted eigenvalues based on \eqref{Lambda} are also shown on the plot.  The curves agree well indicating that the new formulation of PPCA accurately predicts the dependence of the eigenvalues on both $\sigma^2_{w_i}$ and $\sigma^2_\epsilon$. 
 
 \vspace{1cm}

\begin{figure}[H]\vspace{1cm}
\includegraphics[width=5.5cm]{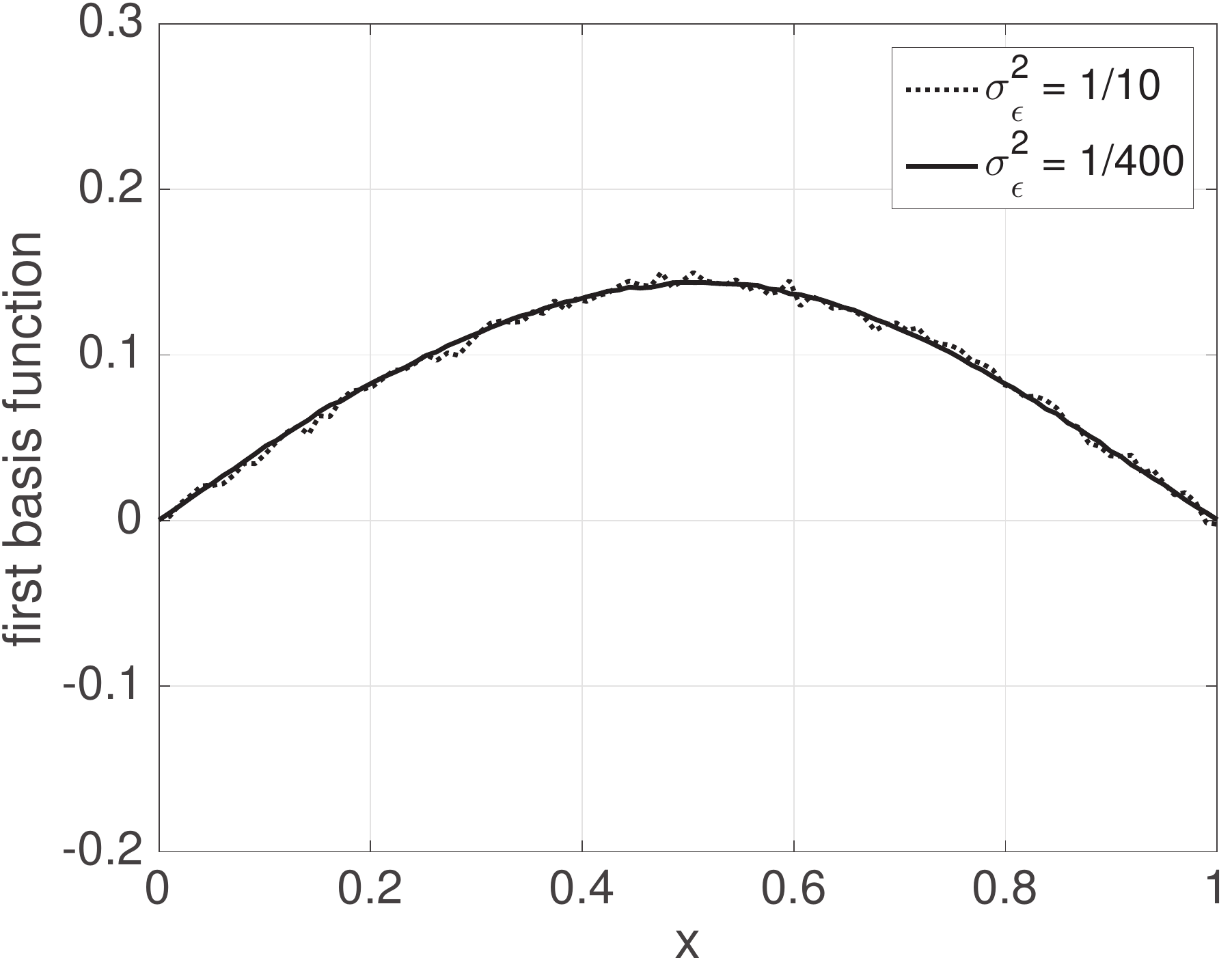}(a)
\includegraphics[width=5.5cm]{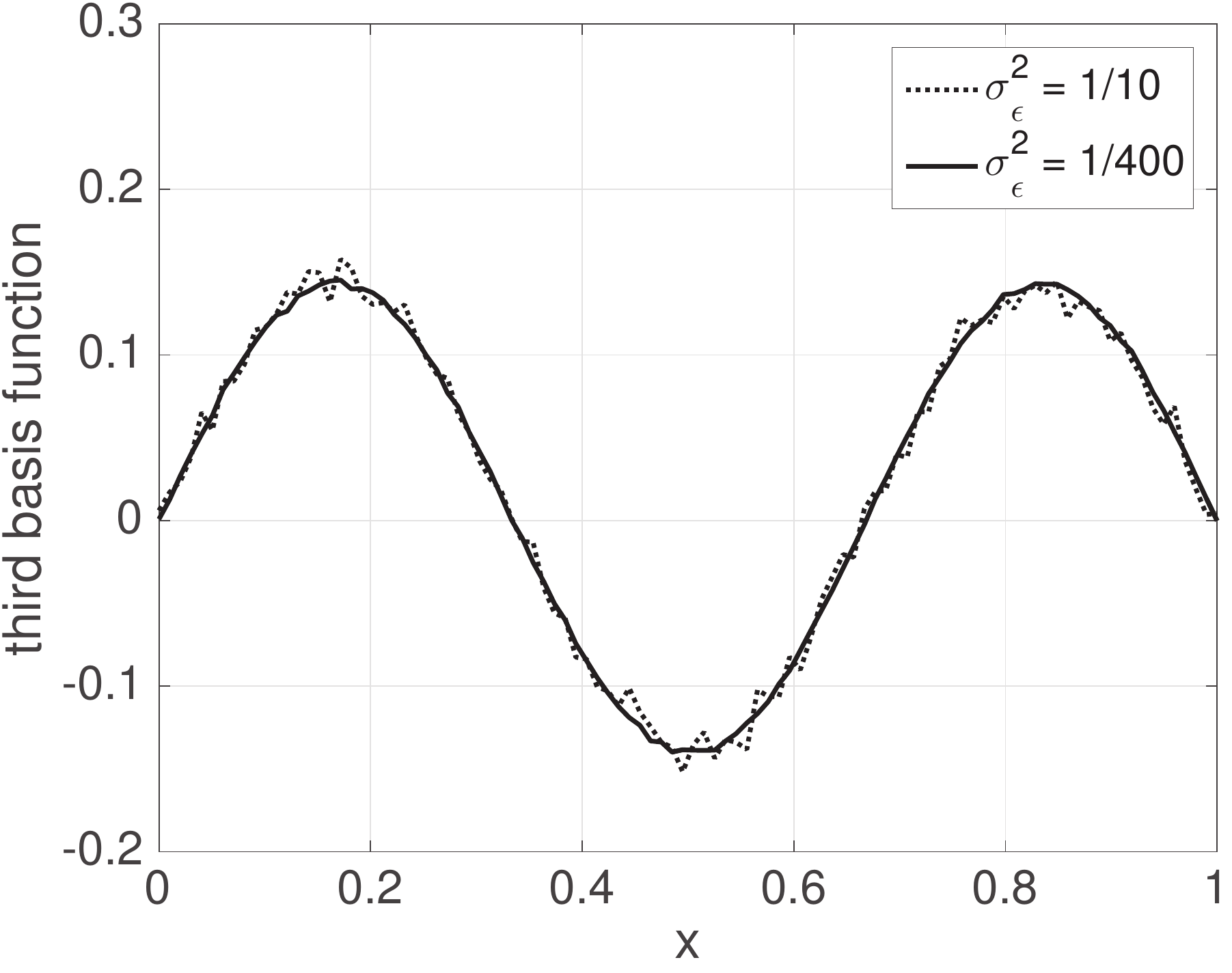}(b)
\includegraphics[width=5.5cm]{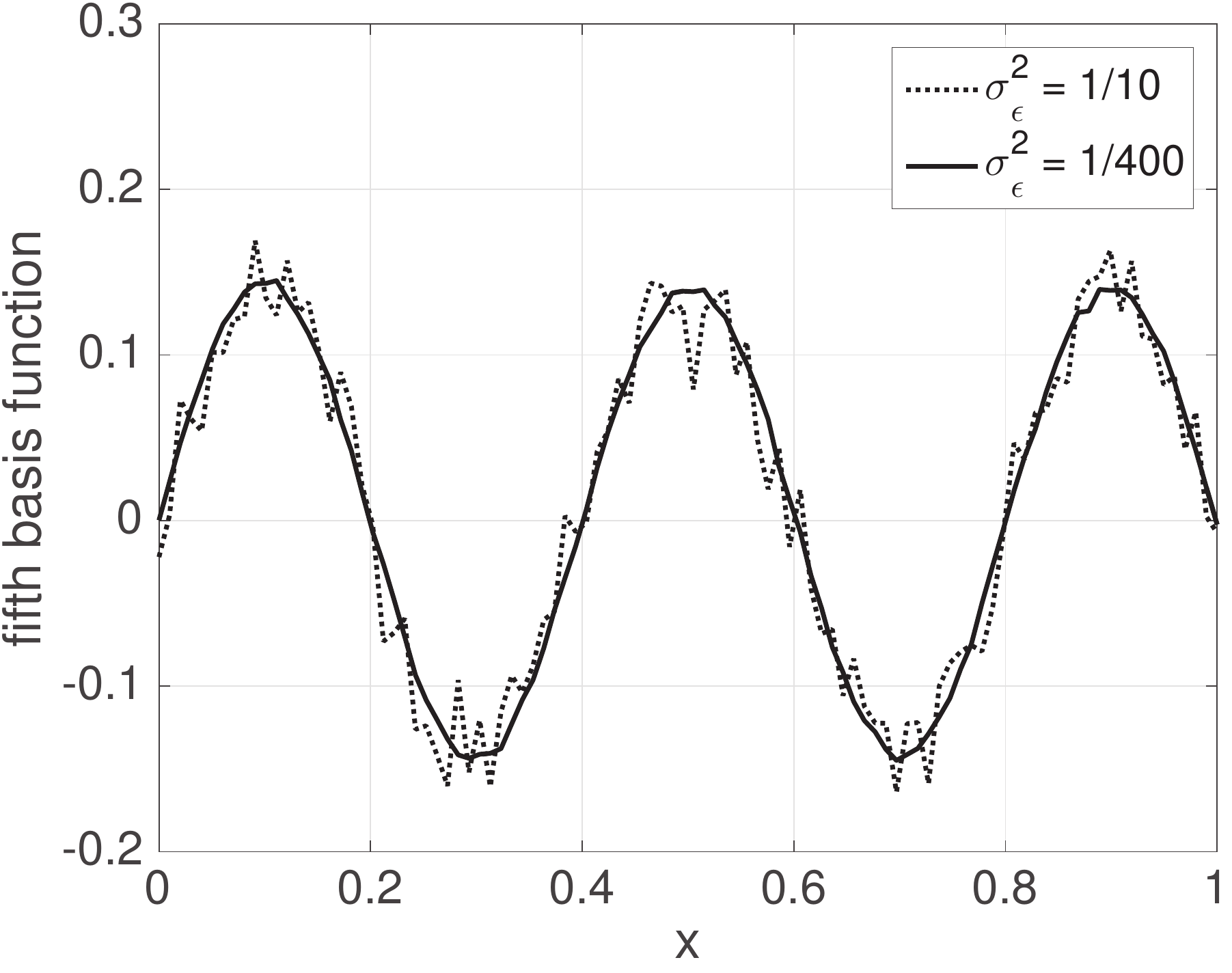}(c)
\includegraphics[width=5.5cm]{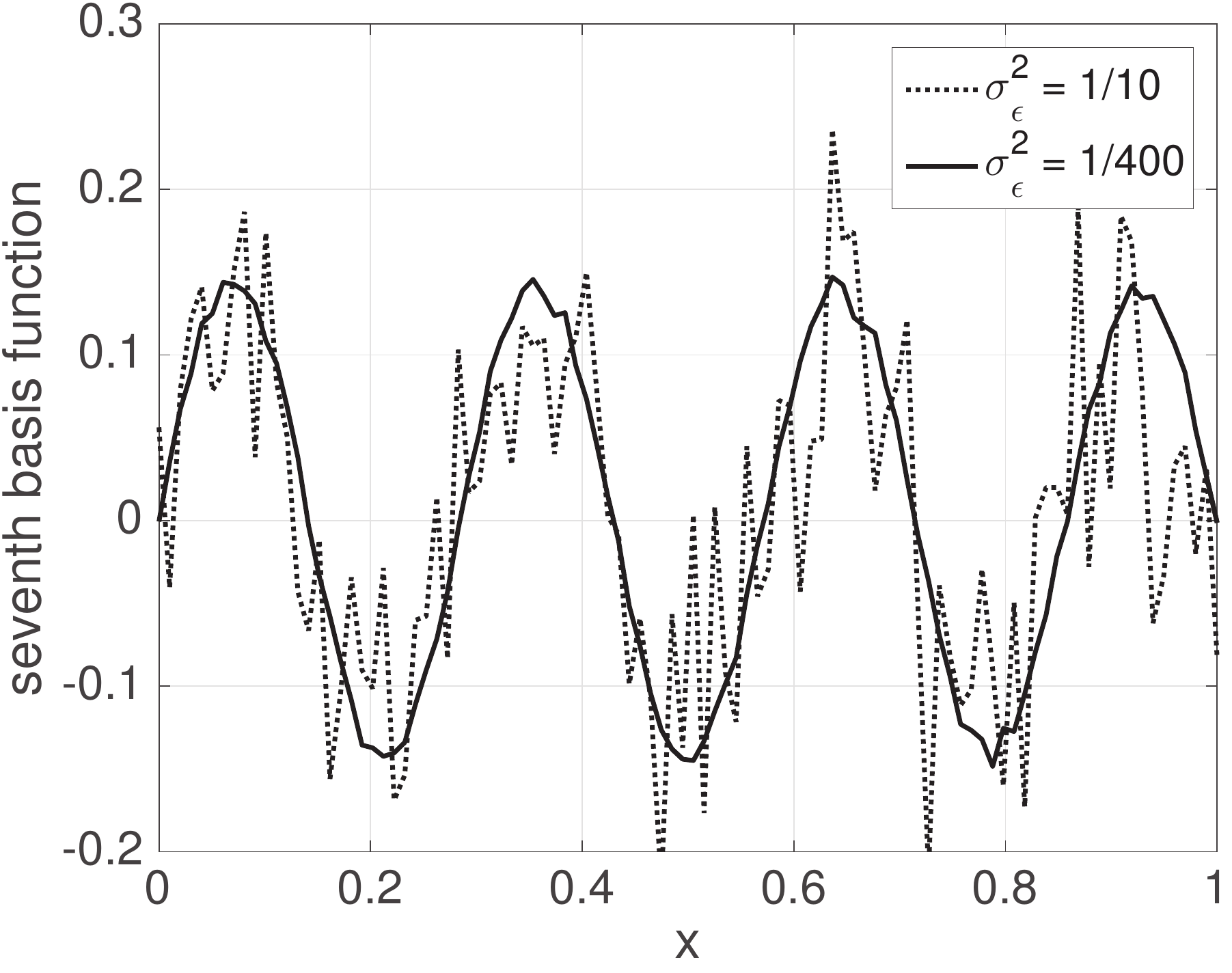}(d)
\caption{Eigenvectors 1,3,5,7 generated using training data with $\sigma_{\epsilon}^{2}$ = 1/10  and $\sigma_{\epsilon}^{2}$ = 1/400.}
\label{basisfun1}
\end{figure}

Fig.~\ref{basisfun1} shows the first, third, fifth, and seventh basis functions for the $\sigma_{\epsilon}^2 = 1/10$ and $\sigma_{\epsilon}^2 = 1/400$ training data. The PPCA is able to extract basis functions that are less affected by noise than the actual data (compare the magnitude of the noise in Fig.~\ref{basisfun1}a for $\sigma_{\epsilon}^2 = 1/10$ with the noise in Fig.~\ref{data_0_1}a).  Even for the 7th mode in the case of $\sigma_\epsilon^2 = 1/10$ shown in Fig.~\ref{basisfun1}d, which had $\sigma_{w_7}^2 = 1/2^6 = 0.016$ being much smaller than $\sigma_\epsilon^2$, the PPCA is able to roughly obtain the correct functional form. This is primarily because of the large number of training data vectors used ($n = 10000$), which allows the PPCA to detect the form of the basis in spite of the numerical noise.  This indicates that one can obtain accurate basis functions for the ROM by either reducing the noise in the data or by increasing the number of data vectors in the training data. 

\vspace{1cm}

\subsection{Model Selection}
\begin{figure}[H]
\includegraphics[width=5.5cm]{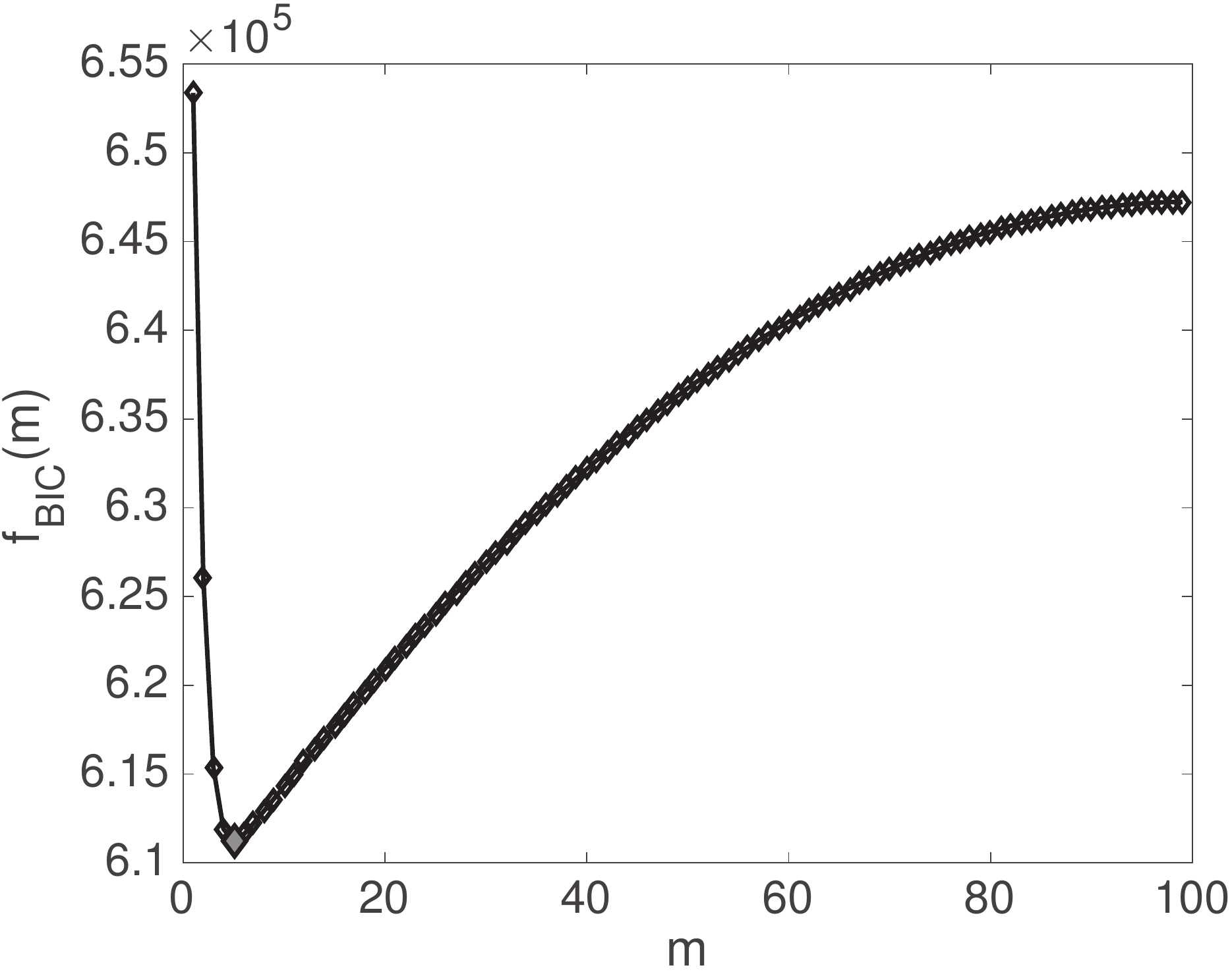}(a)
\includegraphics[width=5.5cm]{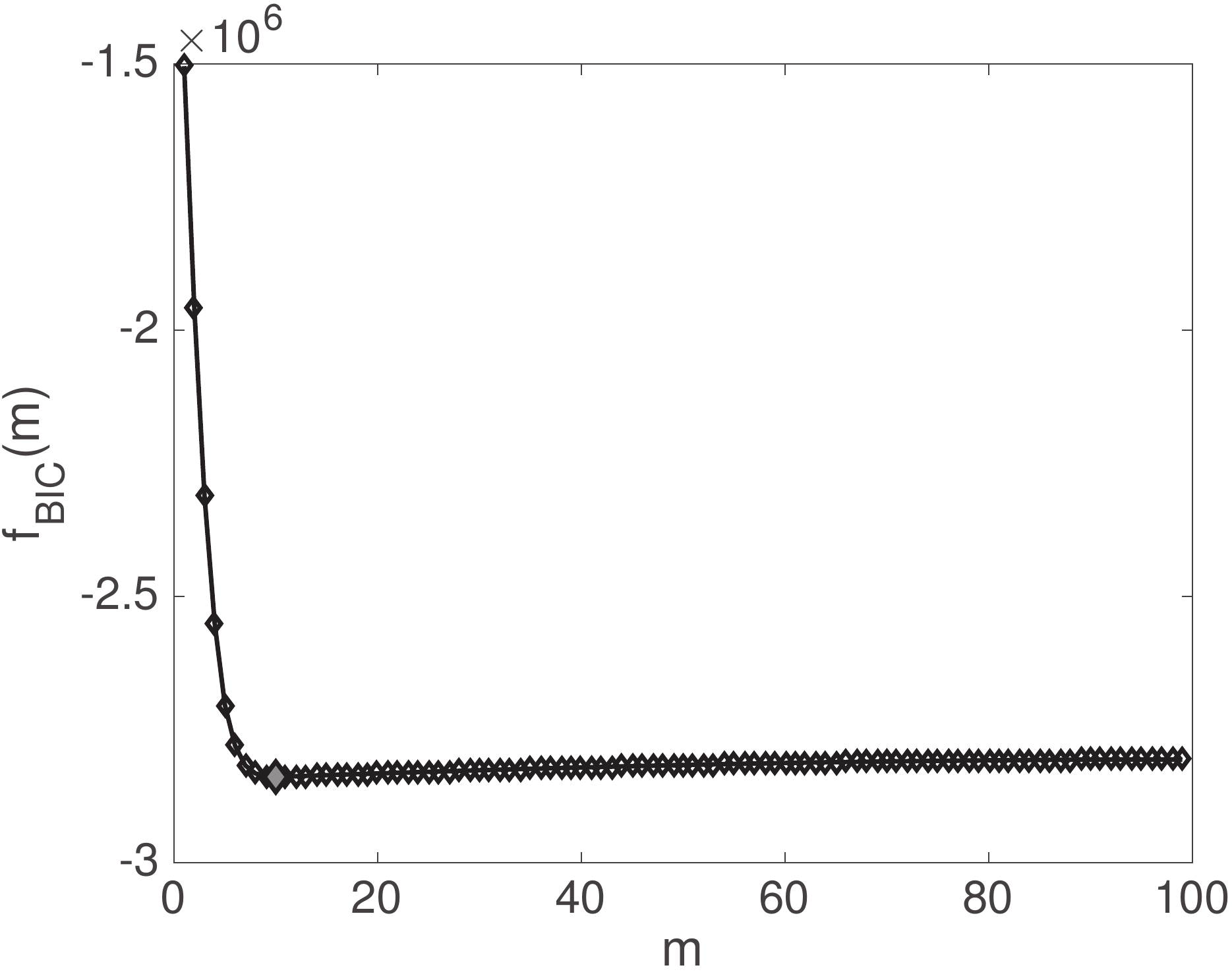}(b)
\caption{BIC as a function of model dimension $m$ (a) for the $\sigma_{\epsilon}^{2} = 1/10$ training data (minimum is at 5 basis functions)  (b) for the $\sigma_{\epsilon}^{2} = 1/400$ training data (minimum is at 10 basis functions).}
\label{BIC_fig}
\end{figure}

As described in Section \ref{dimension}, the Bayesian information criteria (BIC) given in \eqref{BIC} is used determine the dimensionality $m$.  Fig.~\ref{BIC_fig} shows $f_{BIC}(m)$ as a function of $m$ for the two training data sets. The minimum of $f_{BIC}(m)$ occurs at $m=5$ and $m=10$ for the $\sigma_{\epsilon}^2 = 1/10$ and $\sigma_{\epsilon}^2 = 1/400$ training data respectively.   Examining fig.~\ref{eig}, it seems that BIC chooses $m$ at the point of the change in decay rate of the eigenvalues.  For $\sigma_{\epsilon}^2 = 1/10$, this occurs at $m=5$ even though only the first three latent variables have $\sigma_{w_i}^2 > \sigma_{\epsilon}^2$.  For $\sigma_{\epsilon}^2 = 1/400$, this occurs at $m=10$, even though only the first eight latent variables have $\sigma_{w_i} > \sigma_{\epsilon}^2$.

Once the number of latent variables is estimated, \eqref{sig} can be used to estimate $\sigma_{\epsilon}^2$ and \eqref{var_latent} to estimate $\sigma_{w_i}^2$. 

Table~\ref{var_less} summarizes the estimated and true values of $\sigma_{\epsilon}^2 $ and $\sigma_{w_i}^2$ for the two sets of training data.   There are only five values of $\sigma_{w_i}^2$ estimated for the training data with $\sigma_\epsilon^2 = 1/10$ because BIC estimated that $m = 5$ for this data.  The first thing to observe from the data is that the estimates of $\sigma_{w_i}^2$ are accurate;  the percentage errors are less than 5\% in all cases and in most cases the error is less than 2\%.  This is also true of the estimates of $\sigma_\epsilon^2$;  the errors are less than 6\% for both training data sets.  This confirms that the new formulation of PPCA correctly estimates $\sigma_{w_i}^2$.  

Another observation from the estimated values of $\sigma_{w_i}^2$ is that the errors do not vary significantly between the $\sigma_{\epsilon}^2$ = 1/400 and 1/10 data sets.  This implies that for this particular data set, the main source of error is not the noise in the measurements, but rather the number of data vectors used to create the data set ($n$).  Both data sets used $n = 10,000$ so the accuracy of the estimates of $\sigma_{w_i}^2$ are similar.  

\begin{table}
\begin{center}
    \begin{tabular}{ c  c  c  c }
    	 & True  & $\sigma_{\epsilon}^2$ = 1/400 & $\sigma_{\epsilon}^2$ = 1/10 \\ \hline

    $\sigma_{\epsilon}^2$ 	& - 		& 0.002497		& 0.100297	\\
    $\sigma_{w_{1}}^2$  	& 1.0000	& 0.9926 		& 1.0039	\\ 
    $\sigma_{w_{2}}^2$  	& 0.5000	& 0.4911 		& 0.4930	\\ 
    $\sigma_{w_{3}}^2$  	& 0.2500 	& 0.2460 		& 0.2496	\\ 
    $\sigma_{w_{4}}^2$  	& 0.1250 	& 0.1218 		& 0.1278	\\ 
    $\sigma_{w_{5}}^2$  	& 0.0625 	& 0.0645 		& 0.0671	\\ 
  $\sigma_{w_{6}}^2$  	& 0.0313 	& 0.0307		& - 		\\
$\sigma_{w_{7}}^2$  	& 0.0156 	& 0.0157 		& - 		\\
  $\sigma_{w_{8}}^2$  	& 0.0078 	& 0.0079 		& - 		\\
$\sigma_{w_{9}}^2$  	& 0.0039 	& 0.0039 		& - 		\\
 $\sigma_{w_{10}}^2$  	& 0.0019 	& 0.0020 		& - 		\\
    \end{tabular}
    \end{center}
 \caption{Estimated and true values of $\sigma_{w_i}^2$ and $\sigma_\epsilon^2$ for the two training data sets.}
\label{var_less}
\end{table}

\subsection{Bayesian Projection with Gaussian Prior} 
The goal of the projection process is to use the developed model to make an accurate estimate of the true signal, $\Phi \vec{w} + \vec\mu$, given a ``trial'' data vector $\vec{y}$ from an unknown scenario ($\vec{w}$) containing random noise $\vec{\epsilon}$ with unknown magnitude, $\sigma_{\epsilon_T}^2$.  To test the approach described in Section~\ref{projection}, two trial data vectors were generated, one with $\sigma_{\epsilon_{T}}^{2} = 1/5$ and the other with  $\sigma_{\epsilon_{T}}^{2} = 1/200$.  The Bayesian projection approach was then used to estimate $\vec{w}$ and $\sigma_{\epsilon_T}^2$ using \eqref{MAP3} and \eqref{eq_sigma}.   To understand how the training process affects the final results, projections were performed with the 5 and 10 basis function models created above using the $\sigma_\epsilon^2 = 1/10$ and 1/400 training data respectively. One additional ROM was created using training data with $\sigma_\epsilon^2 = 1.5 \times 10^{-12}$ and $m=50$. BIC selected 42 basis functions for this model.

Projections for four cases were performed.  Three cases correspond to trial data with a large noise magnitude ($\sigma_{\epsilon_T}^2 = 1/5$). For these cases reconstructions were done with the ROMs obtained from the $\sigma_\epsilon^2$ = 1/10, 1/400, and $1.5\times 10^{-12}$ training data.  These cases, allowed us to investigate how the quality of the basis functions and the number of basis functions in the ROM affected the projections.  The last case performed a reconstruction using lower noise trial data ($\sigma_{\epsilon_T}^2 = 1/200$) with the ROM obtained from the $\sigma_\epsilon^2$ = 1/400 data.  This case was used to determine the effect of the magnitude of noise in the trial data.
The combination of low noise in the trial data but high noise in the training data is not of practical interest because it is assumed that the training data will always be of the same or higher quality than the trial data.

\begin{figure}[H]\vspace{1cm}
\centering
\includegraphics[width=5.5cm]{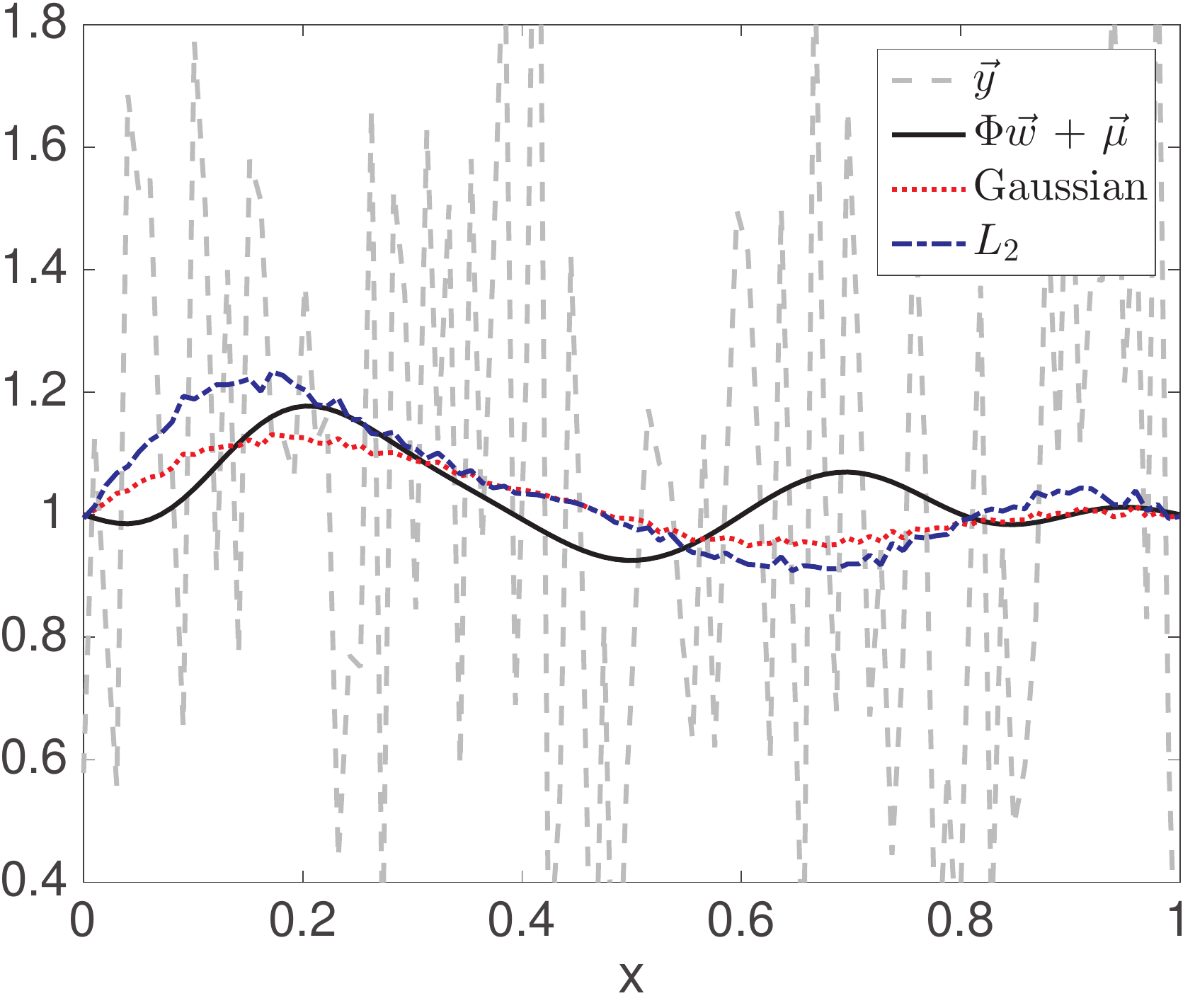}(a)
\includegraphics[width=5.5cm]{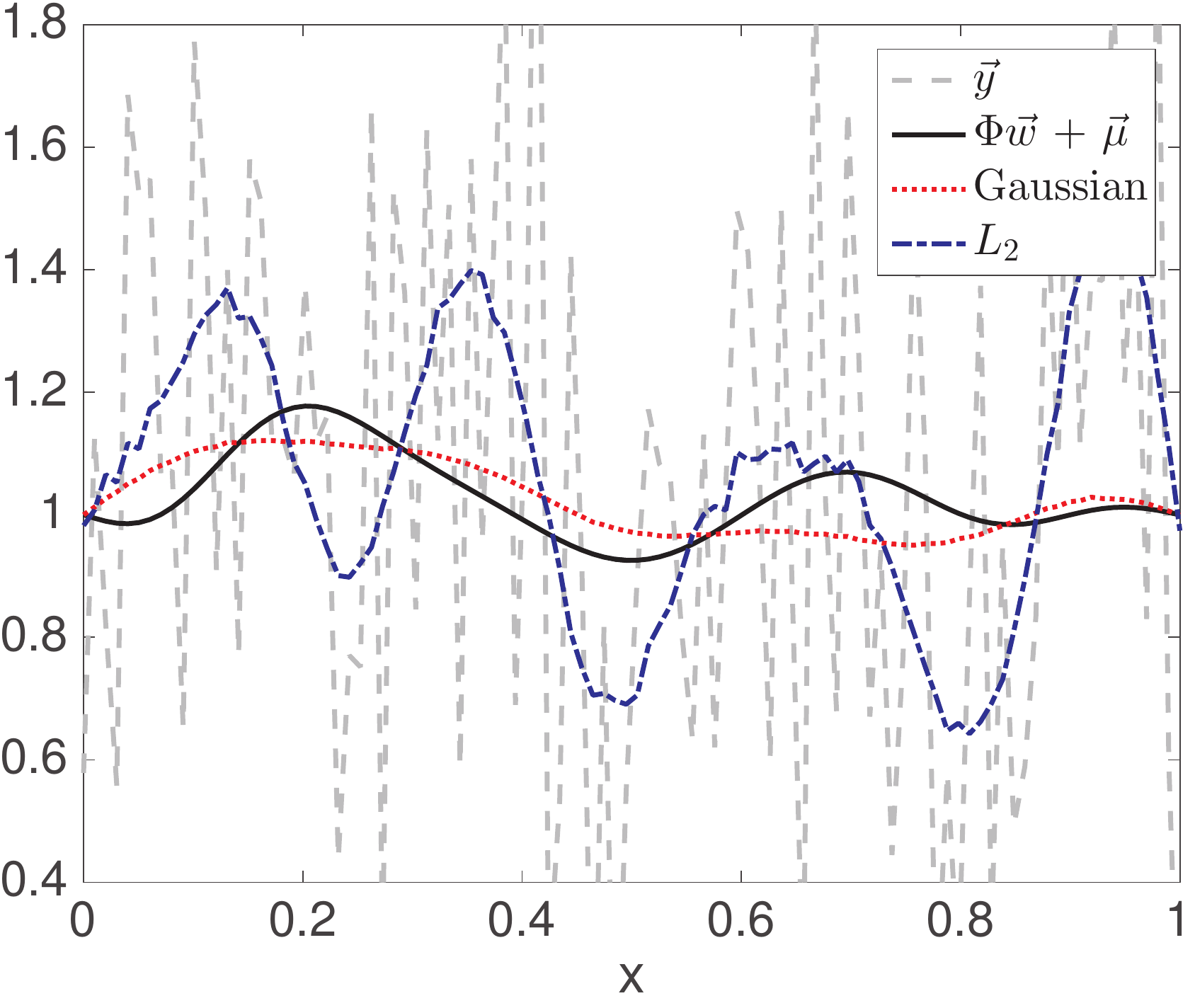}(b) \\
\includegraphics[width=5.5cm]{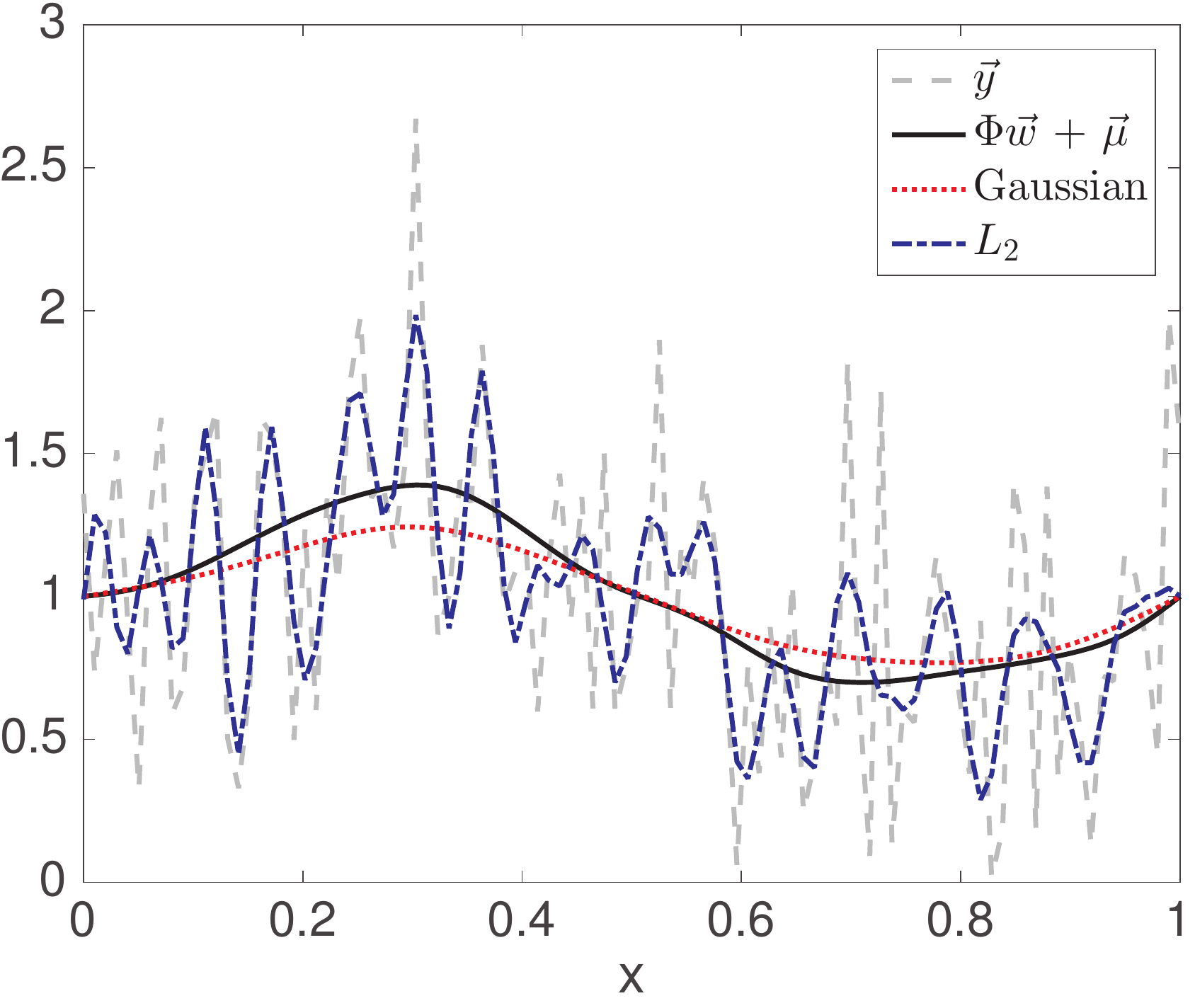}(c)
\includegraphics[width=5.5cm]{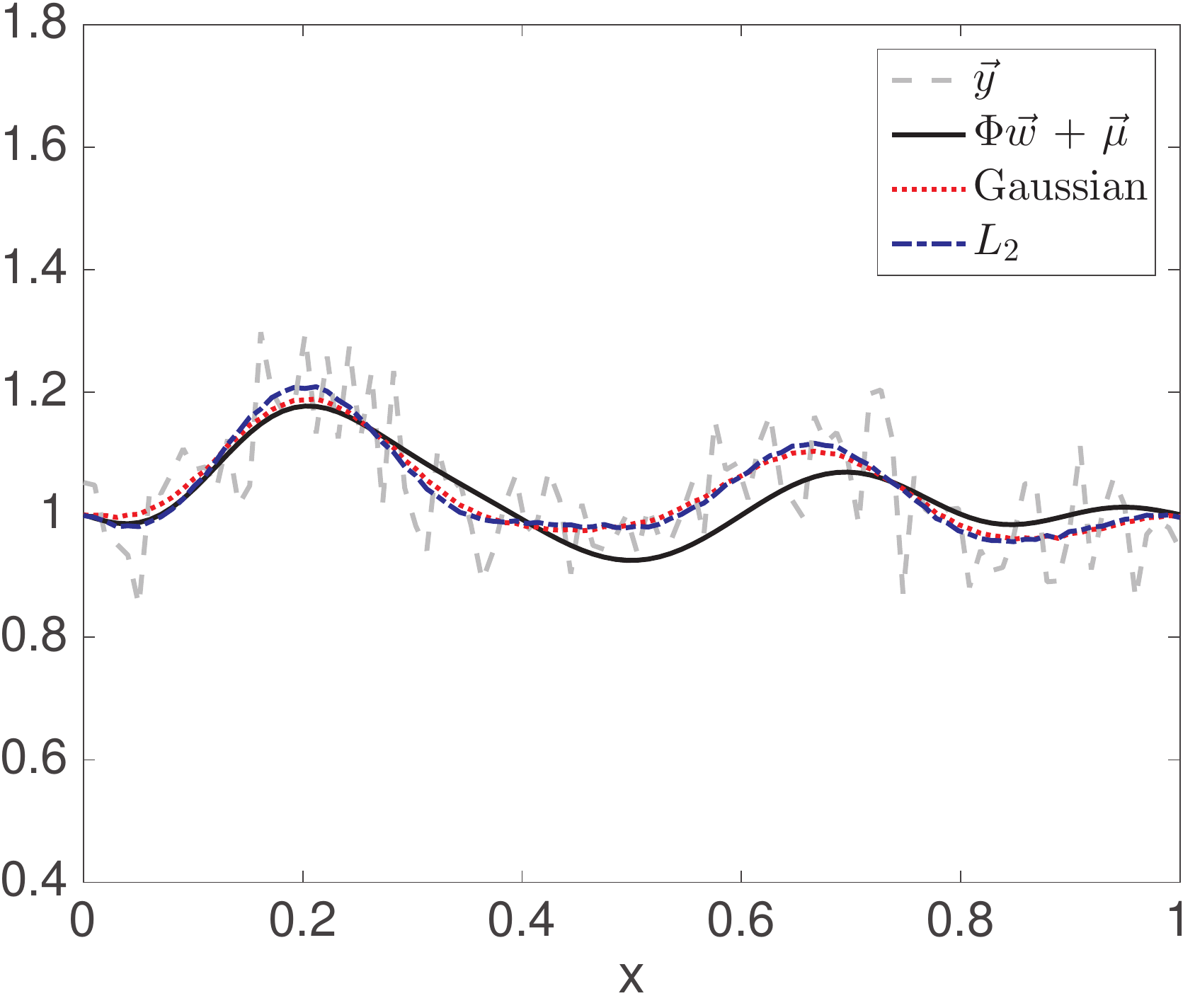}(d)
\caption{Realization $\vec{y}$, true solution $\Phi \vec{w} + \vec{\mu}$, Gaussian-prior projection and $L_2$ projection for the cases: (a) $\sigma_{\epsilon_{T}}^{2} = 1/5$, $\sigma_{\epsilon}^2 = 1/10$, (b) $\sigma_{\epsilon_{T}}^{2} = 1/5$, $\sigma_{\epsilon}^2 = 1/400$,  (c) $\sigma_{\epsilon_{T}}^{2} = 1/5$, $\sigma_{\epsilon}^2 = 1.5 \times 10^{-12}$, (d) $\sigma_{\epsilon_{T}}^{2} = 1/200$, $\sigma_{\epsilon}^2 = 1/400$.}
\label{training_reconstruction}
\end{figure}

Fig.~\ref{training_reconstruction} shows the projection results.  In each figure, the red dotted curve is the ROM projection using a Gaussian prior and the blue dash-dotted curve represents a standard $L_2$ projection. The gray dash line is the trial data, $\vec{y}$, and the solid black line is the true signal i.e. without the added noise.   Fig.~\ref{training_reconstruction}(a) shows the reconstruction of the $\sigma_{\epsilon_T}^2 = 1/5$ trial data using the ROM created using the $\sigma_{\epsilon}^2 = 1/ 10$ training data and Fig.~\ref{training_reconstruction}(b) shows the reconstruction using the ROM created using the $\sigma_{\epsilon}^2 = 1/400$ training data.   Both ROMs significantly reduce the noise in the data, but the $\sigma_{\epsilon}^2 = 1/400$ ROM produces smoother predictions because of the higher quality eigenvectors obtained with the PPCA.  Comparing the $L_2$ projection to the Gaussian-prior projections shows that the Gaussian prior projections give a closer approximation to the true solution.  This is true in both cases, but more so in the $\sigma_{\epsilon}^2 = 1/400$ ROM.  This is because this ROM has more basis functions, which allows the $L_2$ projection to better represent the noise in the data and thus increases the deviation from the true signal.   This is further verified by  Fig.~\ref{training_reconstruction}(c), which shows that with increasing numbers of basis functions in the ROM (42 for this case), the $L_2$ projection result actually deviates from the true solution whereas the Gaussian projection does not.   The projection of the $\sigma_{\epsilon_T}^2 = 1/200$ data shown in Figure~\ref{training_reconstruction}(D) shows that as the noise in the trial data is reduced, both projections approach the noise-free signal, but the Gaussian-prior projection is still more accurate and has fewer spurious oscillation that the $L_2$ projection.  

Table~\ref{var_trial} shows the true and estimated values for $\sigma_{\epsilon_{T}}^{2}$ for the four trial cases shown in Fig.~\ref{training_reconstruction}.   Estimated values are compared for Gaussian and $L_{2}$ projection.   Note that \eqref{eq_sigma} is used to calculate the $\sigma_{\epsilon_{T}}^{2}$ for both projection approaches.  As the number of basis functions in the model increase, the value of $\sigma_{\epsilon_T}^2$ for the Gaussian-prior projection remains relatively constant while the $L_2$ projection approach decreases.  All of the predicted values using the Gaussian-prior projection are within 3\%. 

\begin{table}
\begin{center}
    \begin{tabular}{c  c c c c} 
      True $\sigma_{\epsilon}^2$ & True $\sigma_{\epsilon_{T}}^{2}$ & Predicted $\sigma_{\epsilon_{T}}^{2}$ - Gauss. proj. & Predicted $\sigma_{\epsilon_{T}}^{2}$ -$L_{2}$ proj.\\ \hline
    1/10 & 1/5 & 0.2312 & 0.2294\\ 
 1/400 & 1/5 & 0.2264 & 0.1854\\ 
 $1.5 \times 10^{-12}$& 1/5  & 0.2067& 0.1309\\
 1/400 & 1/200  & 0.00467 &0.00456 \\
    \end{tabular}
\end{center}
\caption{True values and estimated values of the variance of the trial data using Gaussian and $L_{2}$ projection. }
\label{var_trial}
\end{table}

\begin{figure}[H]\vspace{1cm}
\centering
\includegraphics[width=5.5cm]{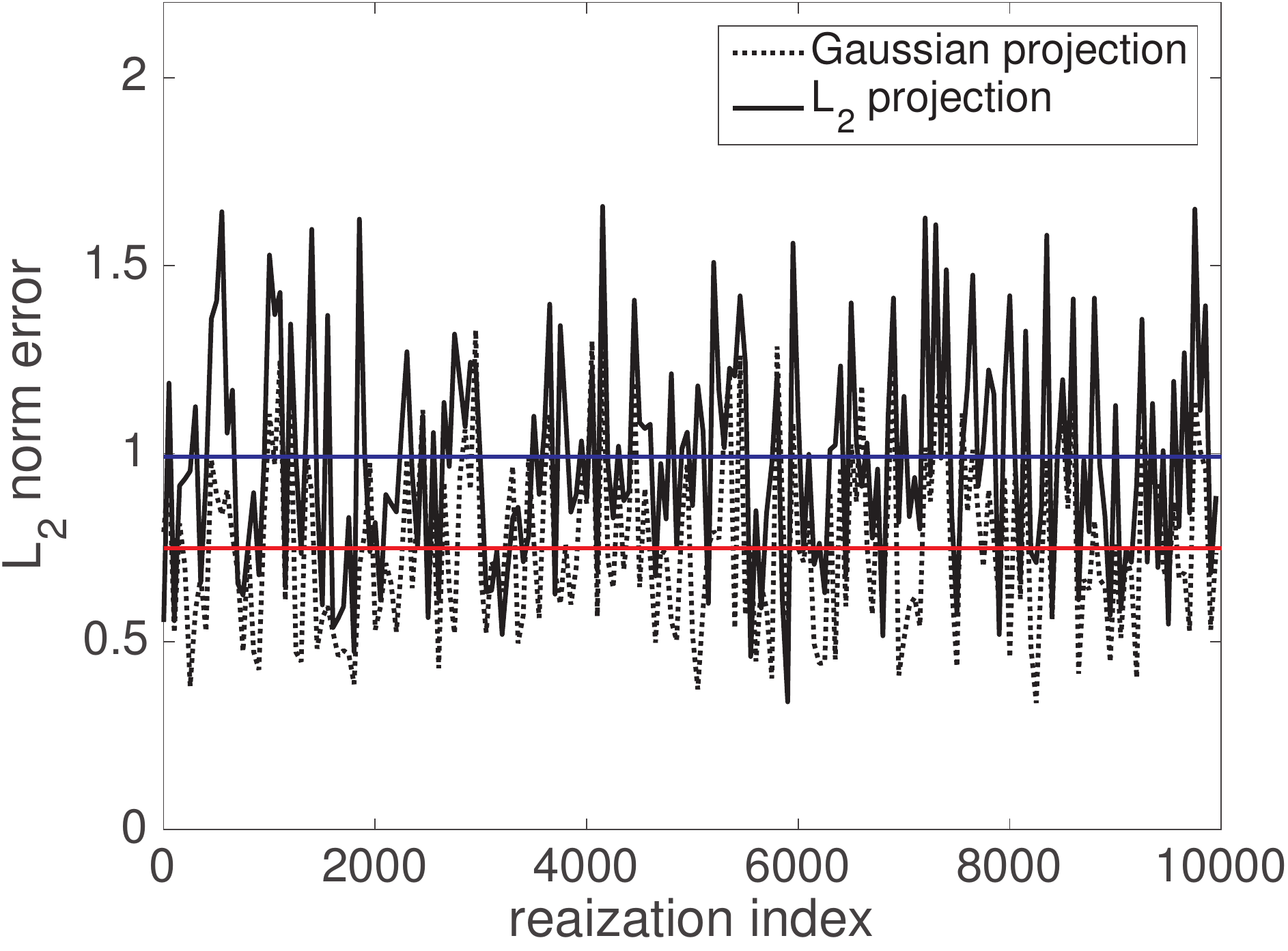}(a)
\includegraphics[width=5.5cm]{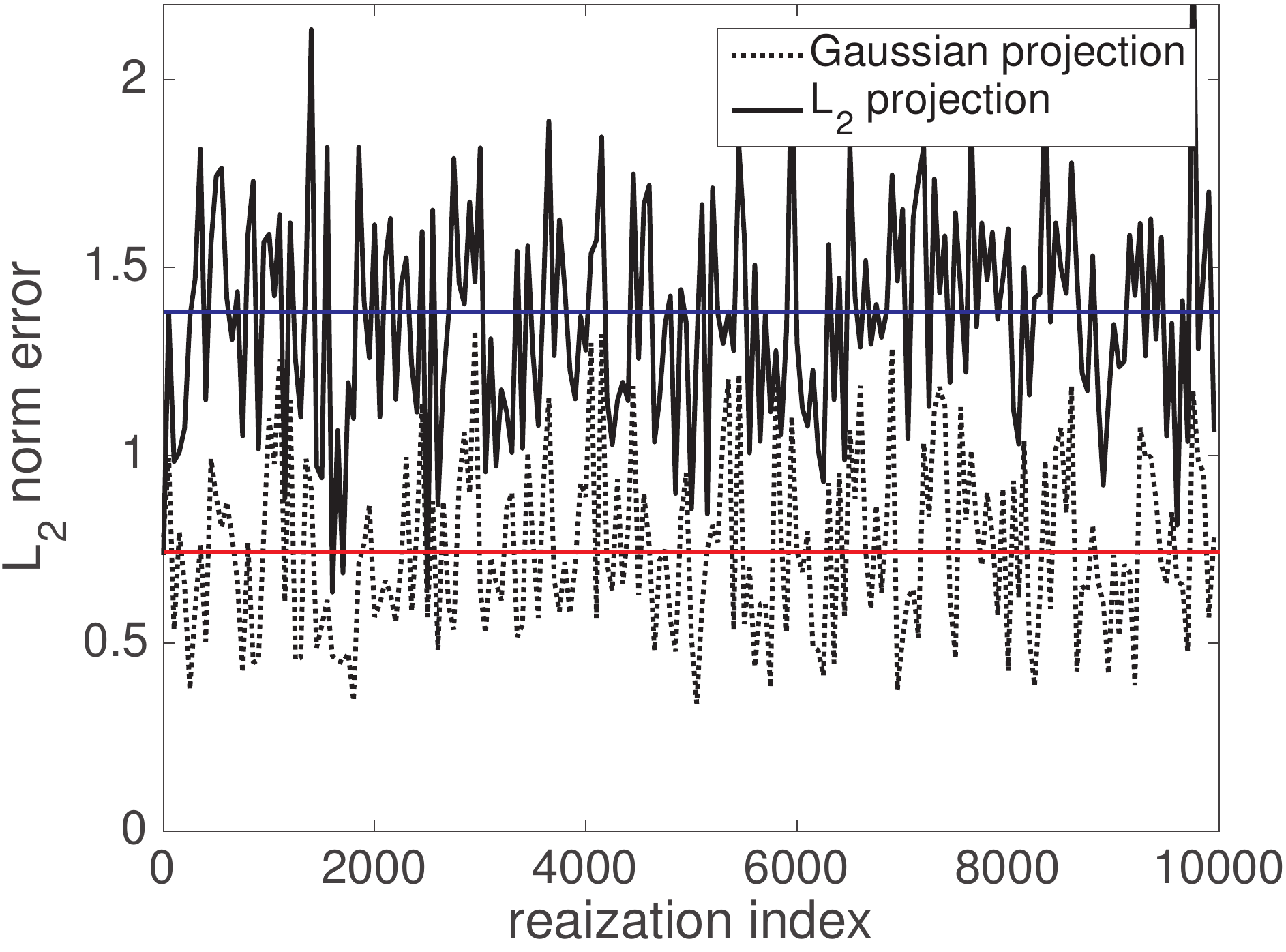}(b) \\
\includegraphics[width=5.5cm]{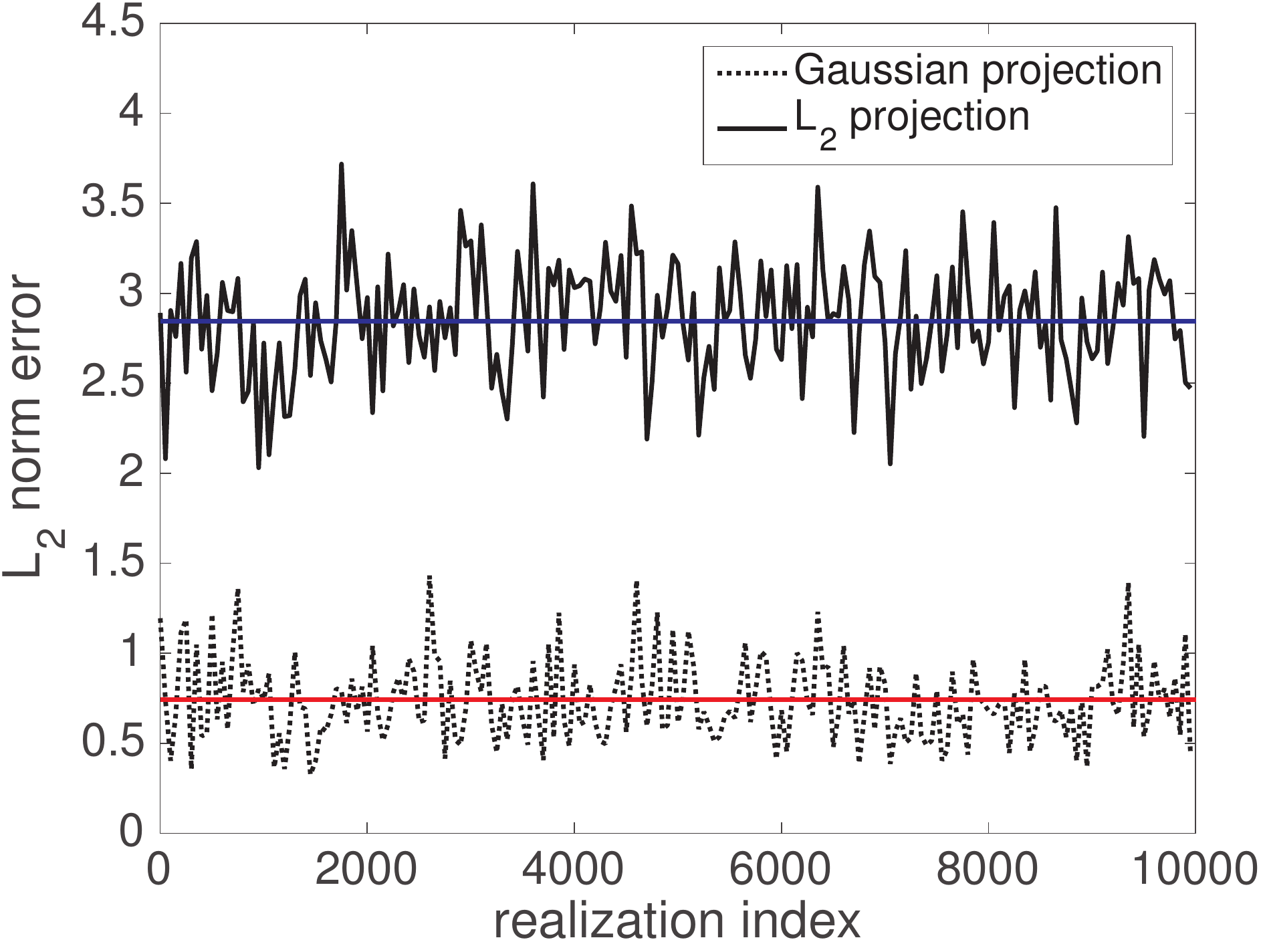}(c)
\includegraphics[width=5.5cm]{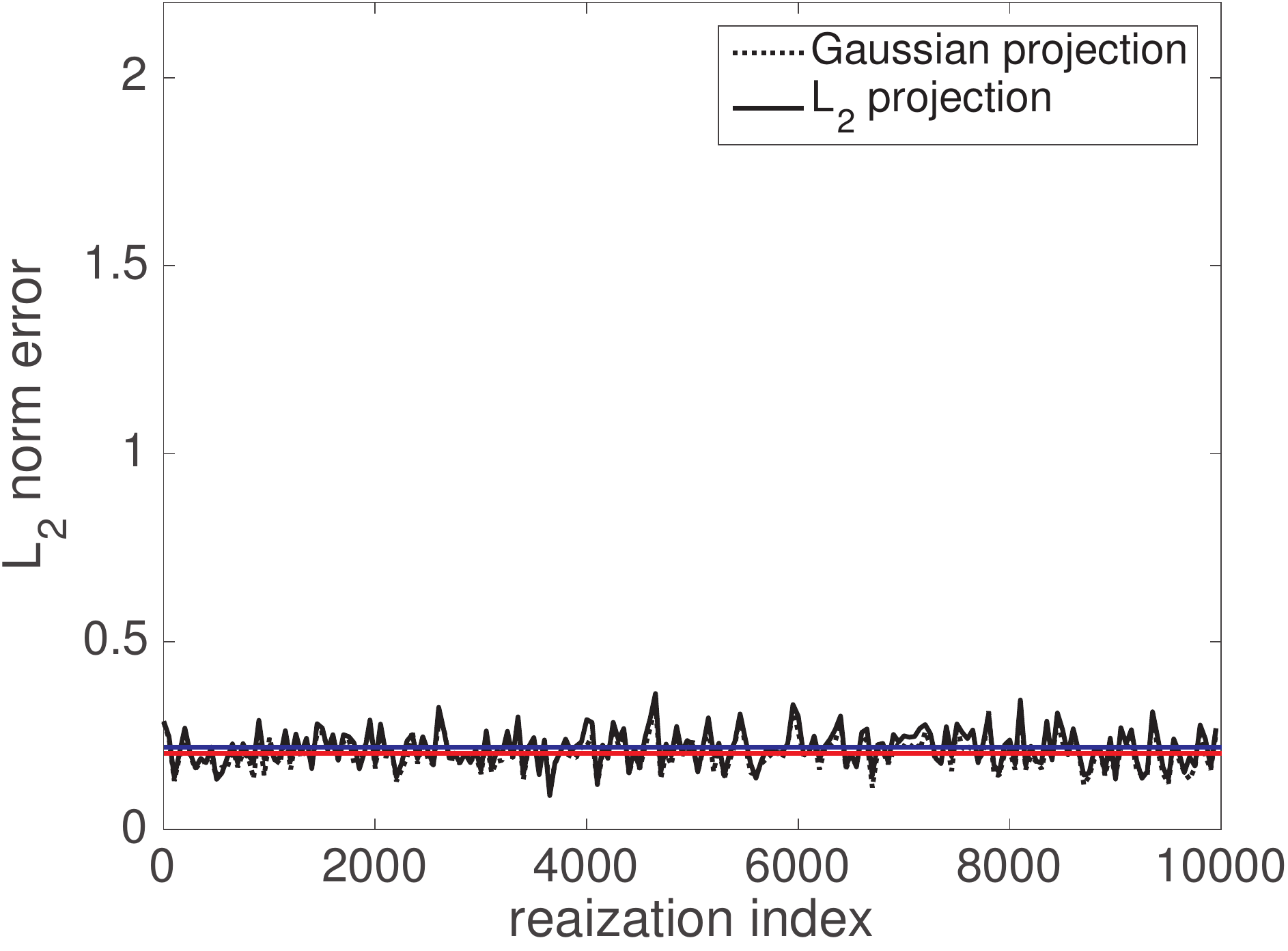}(d)
\caption{$L_{2}$ norm of the error of Gaussian-prior projections and $L_2$ projections for 10000 trial realizations: (a) $\sigma_{\epsilon_{T}}^{2} = 1/5$, $\sigma_{\epsilon}^2 = 1/10$, (b) $\sigma_{\epsilon_{T}}^{2} = 1/5$, $\sigma_{\epsilon}^2 = 1/400$,  (c) $\sigma_{\epsilon_{T}}^{2} = 1/5$, $\sigma_{\epsilon}^2 = 1.5 \times 10^{-12}$, (d) $\sigma_{\epsilon_{T}}^{2} = 1/200$, $\sigma_{\epsilon}^2 = 1/400$. The red line represents the average error for Gaussian-prior projection and the blue line represents the average error for $L_{2}$ projection.}
\label{training_error}
\end{figure}

To verify that the above identified trends are not particular to the trial data vector examined, for 10000 trial data realizations an error was computed by comparing the true signal to the projections.  For the above described cases, the $L_2$ norm of the error,
\begin{equation}
E =  \sqrt{  \left(\vec y_{true} - \vec y_{proj}\right)^T \left(\vec y_{true} - \vec y_{proj}\right) },
\end{equation}
is shown in Fig.~\ref{training_error} for each realization.   The horizontal red and blue curves in the figure are the mean value of the error over the 10000 realizations for the Gaussian-prior projection and the $L_2$ projection respectively.   For all cases, the Gaussian-prior projection on average has less error than the $L_2$ projection. Consistent with the differences seen in Fig.~\ref{training_reconstruction}b and c, this difference is most significant for the case with $\sigma_{\epsilon_{T}}^{2} = 1/5$ using the ROMs with a larger numbers of basis functions. For the case of $\sigma_{\epsilon}^2 = 1/400$, the average error for the Gaussian-prior projection was  0.74 and the average error for $L_2$ projection was  1.38.  For $\sigma_{\epsilon}^2 = 1.5 \times 10^{-12}$, the errors were 0.74 and 2.84.  Note that the scale on Fig.~\ref{training_error}c is larger because the $L_2$ error was large for this case.  Comparing the mean Gaussian-prior projection error between plots a, b, and c shows that the Gaussian-prior projection errors are relatively insensitive to the number of basis functions in the model.  Figure~\ref{training_error}d shows that as $\sigma_{\epsilon_T}^2$ decreases the mean error of both approaches decreases.

\vspace{1cm}

\section{Conclusions}
In this work, statistical approaches were used for each step of the ROM process. Basis functions were generated using the probabilistic formulation of the principal component analysis (PPCA).  The conventional PPCA was improved so that the derived basis functions are orthonormal and the variance of the latent variables is estimated rather than assumed to be one.   These improvements allowed a more intuitive interpretation of the eigenvalues of the PPCA and their relation to the noise in the system and the variance of the latent variables.  

Bayesian model selection was used to select the number of basis functions for the ROM. The Bayesian Information Criteria (BIC) was shown to give reasonable predictions of the number of detectable latent variables in the model, $m$.  Training data with lower noise allowed a larger number of latent variables to be identified, and BIC correspondingly gave a larger estimate for $m$.

Using the information obtained from the improved PPCA to create a Gaussian prior for the latent variables, Bayesian parameter estimation was used  to estimate the latent variables associated with a data vector from a new scenario with an unknown amount of noise.  Equations were derived that simultaneously predict the most-likely latent variable vector and the amount of noise in the new scenario.   The model problem showed that the true (noise-free) signal could be accurately reproduced from noisy data using this approach.   This approach also outperforms  a standard $L_2$ projection to determine the latent variables especially when the noise in the data vector is large and there are many basis functions in the model.

The above formulation provides a comprehensive probabilistic approach for performing reduced order modeling of stochastic systems.  These reduced order models can then be used to provide rapid, accurate predictions for stochastic problems where repeated analyses of similar scenarios must be performed.

\section{Acknowledgements}
This manuscript was authored in part by National Security Technologies, LLC, under contract DE-AC52-06NA25946 with the U.S. Department of Energy and supported by the Site-Directed Research
and Development program. The United States Government retains and the publisher, by accepting
the article for publication, acknowledges that the United States Government retains a nonexclusive,
paid-up, irrevocable, worldwide license to publish or reproduce the published form of this manuscript,
or allow others to do so, for United States Government purposes. The U.S. Department of Energy
will provide public access to these results of federally sponsored research in accordance with the DOE
Public Access Plan (http://energy.gov/downloads/doe-public-access-plan). DOE/NV/25946-{}-3089.

\end{document}